\documentclass[a4paper,12pt]{amsart}
\usepackage{amssymb}
\usepackage{cite}
\usepackage{ifthen}
\usepackage[dvips]{graphicx}
\nonstopmode \numberwithin{equation}{section}
\usepackage{amssymb}
\usepackage{ifthen}
\usepackage{graphicx}
\usepackage{amsmath}
\usepackage[T1]{fontenc} 
\usepackage[utf8]{inputenc}
\usepackage[usenames,dvipsnames]{color}
\usepackage{color}
\usepackage[english]{babel}
\usepackage{fancyhdr}
\usepackage{fancybox}
\usepackage{tikz}

\theoremstyle{plain}
\newtheorem{prop}{Proposition}

\newtheorem{ques}{Question}

\newtheorem{conj}{Conjecture}

\theoremstyle{definition}
\newtheorem{defi}{Definition}[section]

\newtheorem{cor}{Corollary}[section]
\newtheorem{thm}{Theorem}[section]

\newtheorem{lem}{Lemma}[section]
\newtheorem{prob}{Problem}
\newtheorem{rem}{Remark}[section]

\newcounter{minutes}
\divide\time by 60
\newcounter{hours}
\multiply\time by 60
\addtocounter{minutes}{-\time}

\newcounter {own}
\def\theown {\thesection       .\arabic{own}}

\newenvironment{pf}[1][]{%
	\vskip 3mm
	\noindent
	\ifthenelse{\equal{#1}{}}%
	{{\slshape Proof. }}%
	{{\slshape #1.} }%
}%
{\qed\bigskip}

\newcounter{alphabet}

\newcommand{\real}{{\operatorname{Re}\,}}

\def\be{\begin{equation}}
	\def\ee{\end{equation}}

\newcommand{\bee}{\begin{enumerate}}
	\newcommand{\eee}{\end{enumerate}}

\newcommand{\blem}{\begin{lem}}
	\newcommand{\elem}{\end{lem}}
\newcommand{\bthm}{\begin{thm}}
	\newcommand{\ethm}{\end{thm}}
\newcommand{\bcor}{\begin{cor}}
	\newcommand{\ecor}{\end{cor}}
\newcommand{\beg}{\begin{examp}}
	\newcommand{\eeg}{\end{examp}}
\newcommand{\begs}{\begin{examples}}
	\newcommand{\eegs}{\end{examples}}

\newcommand{\bdefn}{\begin{defn}}
	\newcommand{\edefn}{\end{defn}}

\newcommand{\bprob}{\begin{prob}}
	\newcommand{\eprob}{\end{prob}}
\newcommand{\bei}{\begin{itemize}}
	\newcommand{\eei}{\end{itemize}}

\newcommand{\bcon}{\begin{conj}}
	\newcommand{\econ}{\end{conj}}
\newcommand{\bcons}{\begin{conjs}}
	\newcommand{\econs}{\end{conjs}}
\newcommand{\bprop}{\begin{prop}}
	\newcommand{\eprop}{\end{prop}}
\newcommand{\br}{\begin{rem}}
	\newcommand{\er}{\end{rem}}
\newcommand{\brs}{\begin{rems}}
	\newcommand{\ers}{\end{rems}}
\newcommand{\bo}{\begin{obser}}
	\newcommand{\eo}{\end{obser}}
\newcommand{\bos}{\begin{obsers}}
	\newcommand{\eos}{\end{obsers}}
\newcommand{\bpf}{\begin{pf}}
	\newcommand{\epf}{\end{pf}}
\newcommand{\ba}{\begin{array}}
	\newcommand{\ea}{\end{array}}
\newcommand{\beq}{\begin{eqnarray}}
	\newcommand{\beqq}{\begin{eqnarray*}}
		\newcommand{\eeq}{\end{eqnarray}}
	\newcommand{\eeqq}{\end{eqnarray*}}

\begin{document}

\title{Bohr inequalities via proper combinations for a certain class of close-to-convex harmonic mappings}

\author{Molla Basir Ahamed}
\address{Molla Basir Ahamed, Department of Mathematics, Jadavpur University, Kolkata-700032, West Bengal, India.}
\email{mbahamed.math@jadavpuruniversity.in}

\author{Partha Pratim Roy}
\address{Partha Pratim Roy, Department of Mathematics, Jadavpur University, Kolkata-700032, West Bengal, India.}
\email{parthacob2023@gmail.com}

\subjclass[{AMS} Subject Classification:]{Primary 30C45, 30C50, 30C80}

\def\thefootnote{}
\footnotetext{ {\tiny File:~\jobname.tex,
printed: \number\year-\number\month-\number\day,
          \thehours.\ifnum\theminutes<10{0}\fi\theminutes }
}\makeatletter\def\thefootnote{\@arabic\c@footnote}\makeatother
\begin{abstract} 	Let $ \mathcal{H}(\Omega) $ be the class of complex-valued functions harmonic in $ \Omega\subset\mathbb{C} $ and each $f=h+\overline{g}\in \mathcal{H}(\Omega)$, where $ h $ and $ g $ are analytic. In the study of Bohr phenomenon for certain class of harmonic mappings, it is to find a constant $ r_f\in (0, 1) $ such that the inequality 
	\begin{align*}
		M_f(r):=r+\sum_{n=2}^{\infty}\left(|a_n|+|b_n|\right)r^n\leq d\left(f(0), \partial\Omega\right) \;\mbox{for}\;|z|=r\leq r_f,
	\end{align*} 
	where  $ d\left(f(0), \partial\Omega\right) $ is the Euclidean distance between $ f(0) $ and the boundary of $ \Omega:=f(\mathbb{D}) $. The largest such radius $ r_f $ is called the Bohr radius and the inequality $ M_f(r)\leq d\left(f(0), \partial\Omega\right) $ is called the Bohr inequality for the class $ \mathcal{H}(\Omega) $. In this paper, we study Bohr phenomenon for the class of close-to-convex harmonic mappings establishing several inequalities. All the results are proved to be sharp.
\end{abstract}

\maketitle
\pagestyle{myheadings}
\markboth{M. B. Ahamed and P. P. Roy}{Bohr inequalities via proper combinations for certain class of close-to-convex harmonic mappings}

\section{introduction}
Let $ \mathbb{D}:=\{z\in \mathbb{C}:|z|<1\} $ denote the open unit disk in $ \mathbb{C} $. In $1914$ (see \cite{Bohr-1914}), Bohr  established that if $\mathcal{B}$ denote the class of all bounded analytic functions $f(z)=\sum_{n=0}^{\infty}a_nz^n$ on $\mathbb{D}$, then the following inequality holds
\begin{align}\label{Eq-1.1}
	B_0(f,r):=|a_0|+\sum_{n=1}^{\infty}|a_n|r^n\leq||f||_\infty:=\displaystyle\sup_{z\in\mathbb{D}}|f(z)| \;\;\mbox{for}\;\; 0\leq r\leq\frac{1}{3},
\end{align}
where $a_k=f^{(k)}(0)/k!$ for $k\geq 0$. The constant $1/3$ is known the Bohr radius and the inequality in \eqref{Eq-1.1} is called the Bohr inequality for the class $\mathcal{B}$. Henceforth, if there exists a positive real number $r_0$ such that an inequality of the form \eqref{Eq-1.1} holds for every elements of a class $\mathcal{F}$ for $0\leq r\leq r_0$ and fails when $r>r_0$, then we shall say that $r_0$ is an sharp bound for $r$ in the inequality w.r.t. to class $\mathcal{F}$.\vspace{1.2mm}

The origin of Bohr phenomenon lies in this seminal work by Harald Bohr (see \cite{Bohr-1914}) for the analytic functions and this classical result found an application
to the characterization problem of Banach algebras satisfying the von Neumann inequality (see \cite{Dixon & BLMS & 1995}). In recent years the
Bohr inequality has attracted many researchers’ attention to study in different classes of the functions.\vspace{1.5mm}

  In fact, study of Bohr inequality for different classes of functions with various functional settings becomes a subject of great interests during past several years and an extensive research work has been done by many authors (see e.g., \cite{Alkhaleefah-PAMS-2019,Paulsen-PLMS-2002,Paulsen-PAMS-2004,Paulsen-BLMS-2006,Ismagilov-2020-JMAA,Kayumov-CRACAD-2018,Kayumov-Khammatova-JMAA-2021,Kayumov-MJM-2022,Kay & Pon & AASFM & 2019,Lata-Singh-PAMS-2022,Liu-JMAA-2021,Liu-Ponnusamy-MN,Liu-Ponnusamy-PAMS-2021,Ponnusamy-HJM-2021} and references therein). The Bohr phenomenon for Hardy space both in single and several variables-along with some Schwarz-Pick
type estimates are established in \cite{Bene-2004}. The Bohr inequality for holomorphic mappings with a lacunary series in several complex variables are established recently in \cite{Lin-Liu-Ponn-AMS-2023}. For different aspects of Bohr phenomenon including multidimensional Bohr inequality, the readers are referred to the articles \cite{Aizn-PAMS-1999,Aizenberg-SM-2005,Aizenberg-AMP-2012,Allu-CMB-2022,Allu-JMAA-2021,Kumar-PAMS-2022,Ponnusamy-JMAA-2022,Bhowmik-2018,Boas-Khavinson-1997,Dixon & BLMS & 1995,Defant-AM-2012,Evdoridis-IM-2018,Galicer-TAMS-2021} and references therein. However, Bohr phenomenon for the classes of harmonic mappings was initiated first in \cite{Abu-CVEE-2010} and was investigated
in \cite{Kay & Pon & Sha & MN & 2018} and subsequently by a number of authors in recent years in \cite{Aha-CMFT-2022,Aha-Allu-BMMSS-2022,Ahamed-AMP-2021,Ahamed-CVEE-2021,Huang-Liu-Ponnusamy-MJM-2021,Kayumov-2018-JMAA,Liu-Ponnusamy-BMMS-2019,Liu-Ponn-Wang-RACSAM-2020}. The recent survey articles \cite{Ali-Abu-Ponn-2016,Ponnusamy-Vijayakumar} and references therein may be good sources for this topic.

\subsection{Rogosinski radius.} Just like the Bohr radius, there exists a concept known as the Rogosinski radius \cite{Rogosinski-1923}, defined as follows: if $f(z)=\sum_{n=0}^{\infty}a_nz^n\in\mathcal{B}$, then $|S_M(z)|=|\sum_{n=0}^{M-1}a_nz^n|<1$ for  $|z|<1/2,$, where $1/2$ is best possible. The radius $1/2$ is called the Rogosinski Radius for the family $\mathcal{B}.$ In \cite{Kayumov-Khammatova-JMAA-2021}, Kayumov \emph{et al.} considered the Bohr-Rogosinski inequality
\begin{align}\label{Eq-1.2}
	|f(z)|+\sum_{n=N}^{\infty}|a_n|r^n\leq1=|f(0)|+d(f(0),\partial \mathbb{D})
\end{align}
for $|z|=r\leq R_N,$ where $R_N$ is the positive root of the equation $2(1+r)r^N-(1-r)^2=0$ and $d(f(0),\partial \mathbb{D})$ denotes the Euclidean distance from $f(0)$ to the boundary of the domain $\partial\mathbb{D}.$ If we replace $|f(z)|$ by $|f(0)|$ and $N=1$ in \eqref{Eq-1.2}, then we see the connection between the Bohr-Rogosinski inequality and classical Bohr inequality. Recently,  Kumar and Sahoo\cite{Kumar-Sahoo-Lova-2022} have generalized the Bohr-Rogosinski inequality \eqref{Eq-1.2} proving them sharp.\vspace{1.5mm}

The objective of the paper is to study Bohr phenomenon with suitable settings in order to establish certain harmonic analogue of some Bohr inequality valid for analytic functions on unit disk.
\begin{defi}\cite{Grigoryan-POTENTIAL-2023}
	Let $A\subset\hat{\mathbb{C}}$ be a set, a function $f: A\rightarrow\mathbb{C}$ is said to be harmonic in an open set $\Omega$ provided $\Omega\subset A,$ $f$ is continuous in $\Omega$ and $f$ is twice continuously differentiable in $\Omega\setminus\{\infty\}$ and satisfies the Laplace equation 
	\begin{align*}
		\dfrac{\partial^2{f(z)}}{\partial x^2}+\dfrac{\partial^2{f(z)}}{\partial y^2}=0,\; z=x+iy\in \Omega\setminus\{\infty\}.
	\end{align*}
	In particular, $f$ is said to be harmonic provided $A$ is an open set and $f$ is harmonic in $A$.
\end{defi}
 Harmonic mappings are instrumental because of its applications in many Science and Engineering branches. To be specific, methods of harmonic mappings have been applied to study and solve the fluid flow problems (see \cite{Aleman-2012,Constantin-2017}). For example, in 2012, Aleman and Constantin \cite{Aleman-2012} established a connection between harmonic mappings and ideal fluid flows. In fact, Aleman and Constantin have developed an ingenious technique to solve the incompressible two dimensional Euler equations in terms of univalent harmonic mappings (see \cite{Constantin-2017} for details).
\subsection{Bohr inequality for harmonic mappings}
 Let $ \mathcal{H}(\Omega) $ be the class of complex-valued functions harmonic in $ \Omega $. It is well-known that functions $ f $ in the class $ \mathcal{H}(\Omega) $ has the following representation $ f=h+\overline{g} $, where $ h $ and $ g $ both are analytic functions in $ \Omega $. The famous Lewy's theorem (see \cite{Lew-BAMS-1936}) states that a harmonic mapping $ f=h+\overline{g} $ is locally univalent on $ \Omega $ if, and only if, the determinant $ |J_f(z)| $ of its Jacobian matrix $ J_f(z) $ does not vanish on $ \Omega $, where
\begin{equation*}
	|J_f(z)|:=|f_{z}(z)|^2-|f_{\bar{z}}(z)|^2=|h^{\prime}(z)|^2-|g^{\prime}(z)|^2\neq 0.
\end{equation*}
In view of this result, a locally univalent harmonic mapping  is sense-preserving if $ |J_f(z)|>0 $ and sense-reversing if $|J_{f}(z)|<0$ in $\Omega$. For detailed information about the harmonic mappings, we refer the reader to \cite{Clunie-AASF-1984,Duren-2004}. In \cite{Kay & Pon & Sha & MN & 2018}, Kayumov \emph{et al.} first established the harmonic extension of the classical Bohr theorem, since then investigating on the Bohr-type inequalities for certain class of harmonic mappings becomes an interesting topic of research in geometric function theory.\vspace{1.2mm}
 
  
We see that \eqref{Eq-1.1} can be written as
\begin{equation}\label{e-1.2}
	d\left(\sum_{n=0}^{\infty}|a_nz^n|,|a_0|\right)=\sum_{n=1}^{\infty}|a_nz^n|\leq 1-|f(0)|=d(f(0),\partial (\mathbb{D})).
\end{equation}
More generally, a class $ \mathcal{F} $ of analytic functions $ f(z)=\sum_{n=0}^{\infty}a_nz^n $ mapping $ \mathbb{D} $ into a domain $ \Omega $ is said to satisfy a Bohr phenomenon if an inequality of type \eqref{e-1.2} holds uniformly in $ |z|\leq \rho_0 $, where $ 0<\rho_0\leq 1 $ for functions in the class $ \mathcal{F} $. Similar definition makes sense for classes of harmonic mappings (see \cite{Kay & Pon & Sha & MN & 2018}) also.\vspace{1.2mm}

In view of the distance formulation of the Bohr inequality as mentioned in \eqref{e-1.2}, Abu-Muhanna (see \cite{Abu-CVEE-2010}) have established the following result for subordination class $ S(f) $ in case of when $ f $ is univalent function showing that the corresponding inequality is sharp.
\begin{thm}\cite{Abu-CVEE-2010}
	If $ g(z)=\sum_{n=0}^{\infty}b_nz^n\in S(f) $ and $ f(z)=\sum_{n=0}^{\infty}a_nz^n $ is univalent, then 
	\begin{align*}
		\sum_{n=0}^{\infty}|b_nz^n|\leq d(f(0), \partial\Omega)\; \mbox{for}\; |z|=\rho_0^*\leq 3-\sqrt{8}\approx 0.17157,
	\end{align*}
 where $ \rho_0^* $ is sharp for Koebe function $ f(z)=z/(1-z)^2 $.
\end{thm}

\par Let $ \mathcal{H} $ be the class of all complex-valued harmonic functions $ f=h+\bar{g} $ defined on the unit disk $ \mathbb{D} $, where $ h $ and $ g $ are analytic in $ \mathbb{D} $ with the normalization $ h(0)=h^{\prime}(0)-1=0 $ and $ g(0)=0 $. Let $ \mathcal{H}_0 $ be defined by $ 	\mathcal{H}_0=\{f=h+\bar{g}\in\mathcal{H} : g^{\prime}(0)=0\}. $ Therefore, each $f=h+\overline{g}\in \mathcal{H}_{0}$ has the following representation 
\begin{equation}\label{e-1.3}
	f(z)=h(z)+\overline{g(z)}=\sum_{n=1}^{\infty}a_nz^n+\overline{\sum_{n=1}^{\infty}b_nz^n}=z+\sum_{n=2}^{\infty}a_nz^n+\overline{\sum_{n=2}^{\infty}b_nz^n},
\end{equation}
where $ a_1 = 1 $ and $ b_1 = 0 $, since $ a_1 $ and $ b_1 $ have been appeared in later results and corresponding proofs.\vspace{1.2mm}

For a class of harmonic mappings, the Bohr radius and its corresponding Bohr inequality analogue in view of distance formulation, similar to that for analytic functions, are defined as follows.
\begin{defi}
Let $ f\in\mathcal{H}_0 $ be given by \eqref{e-1.3}. In the study of it is the Bohr phenomenon is to find a constant $ r_f\in (0, 1] $ such that the inequality 
	\begin{align*}
		r+\sum_{n=2}^{\infty}\left(|a_n|+|b_n|\right)r^n\leq d\left(f(0), \partial\Omega\right) \;\mbox{for}\;|z|=r\leq r_f.
	\end{align*} 
The largest such radius $ r_f $ is called the Bohr radius for the class $ \mathcal{H}_0 $.
\end{defi}
	Before stating the main result of the paper, we need some background of Bohr phenomenon for the class $\mathcal{B}$ of all analytic self-maps on $\mathbb{D}.$ In the section, we discuss in details several refined and improved Bohr inequalities established in recent years.
\section{Bohr inequalities for the class $\mathcal{B}$}
To state our main results, we need some preparation. We fix some notations here 
\begin{align*}
	\mathcal{B}=\{f(z)=\sum_{n=0}^{\infty}|a_n|z^n\in \mathcal{B} \;:\;||f||_{\infty}\leq 1\}.
\end{align*} 
For $f\in \mathcal{B}$ and $f_0(z):=f(z)-f(0)$, we let for convenience 
\begin{align*}
	B_k(f,r):=\sum_{n=k}^{\infty}|a_n|r^n \;\mbox{for }k\geq 0,\; ||f_0||^2_r:=\sum_{n=1}^{\infty}|a_n|^2r^{2n}.
\end{align*}
Further, for $ f\in\mathcal{B} $, we define the quantity $A(f_0,r)$ by 
\begin{align*}
		A(f_0,r):=\left(\dfrac{1}{1+|a_0|}+\dfrac{r}{1-r}\right)||f_0||^2_r.
\end{align*}
With reference to Rogosinski's inequality and radius investigated in \cite{Rogosinski-1923}, Kayumov and Ponnusamy \cite{Kayumov-Khammatova-JMAA-2021} have introduced and derived the Bohr-Rogosinski inequalities and Bohr-Rogosinski radius for the class $ \mathcal{B} $.
\begin{thm}\cite{Kayumov-Khammatova-JMAA-2021}\label{Th-2.1}
	Suppose that $ f(z)=\sum_{n=0}^{\infty}a_nz^n\in\mathcal{B} $. For $N,m\in\mathbb{N}, $ 
	\begin{align*}
		|f(z^m)|+B_N(f,r)\leq 1\;\;\mbox{for}\;\; r\leq R_{m,N}, 
	\end{align*} 
	where $ R_{m,N} $ is the positive root of the equation $ 2(1+r^m)r^N-(1-r)(1-r^m)=0 $. The number $ R_{m,N} $ cannot be improved. Moreover, 
	$ \lim\limits_{N\rightarrow\infty} R_{m,N}=1$ and 	$ \lim\limits_{m\rightarrow\infty} R_{m,N}=A_N,$ where $A_N$ is the positive root of the equation $2r^N=1-r.$  In addition,  
	for $N,m\in\mathbb{N}, $ 
	\begin{align*}
		|f(z^m)|^2+B_N(f,r)\leq 1\;\;\mbox{for}\;\; r\leq R^{\prime}_{m,N},
	\end{align*} 
	where $ R^{\prime}_{m,N} $ is the positive root of the equation $ (1+r^m)r^N-(1-r)(1-r^m)=0 $. The number $ R^{\prime}_{m,N} $ cannot be improved. Moreover, 
	$ \lim\limits_{N\rightarrow\infty} R^{\prime}_{m,N}=1$ and 	$ \lim\limits_{m\rightarrow\infty} R^{\prime}_{m,N}=A^{\prime}_N,$ where $A^{\prime}_N$ is the positive root of the equation $r^N=1-r.$ 
\end{thm}
\subsection{Refined versions of the Bohr inequality for the class $\mathcal{B}$}
For an extension of the results discussed above, we refer to the recent article by Ponnusamy and Vijayakumar \cite{Ponnusamy-Vijayakumar}. In comparison of $ \sum_{n=0}^{\infty}|a_n|r^n $ with another functional often considered in function theory, namely $ \sum_{n=0}^{\infty}|a_n|^2r^{2n} $ which is abbreviated as $ ||f||^2_r $. As refinement of the classical Bohr inequality, for $ p=1, 2 $, we define
\begin{align*}
	\mathcal{A}_{p, f}(r):= |a_0|^p+B_1(f,r)+A(f_0,r),
\end{align*}
where $ f_0(z):=f(z)-a_0 $. Recently, Ponnusamy \emph{et al.} \cite{Ponnusamy-RM-2020} have obtained the following result as a refinement of the classical Bohr inequality.
\begin{thm}\cite{Ponnusamy-RM-2020}\label{th-11.33}
Suppose that $ f\in\mathcal{B} $ with $ f(z)=\sum_{n=0}^{\infty}a_nz^n $ and $ f_0(z):=f(z)-a_0 $. Then $ \mathcal{A}_{1, f}(r)\leq 1 $ and the numbers $ 1/(1+|a_0|) $ and $ 1/(2+|a_0|) $ cannot be improved. Further, $ \mathcal{A}_{2, f}(r)\leq 1 $ and the numbers $ 1/(1+|a_0|) $ and $ 1/2 $ cannot be improved.
\end{thm}
In what follows, $ \lfloor x \rfloor $ denotes the largest integer no more than $ x $, where $ x  $ is a real number. Recently, Liu \emph{et al.} \cite{Liu-Liu-Ponnusamy-BSM-2021} have obtained the following refined version of the Bohr-Rogosinski inequality also.
\begin{thm}\cite{Liu-Liu-Ponnusamy-BSM-2021}\label{TH-2.3}
	Suppose that $ f\in\mathcal{B} $ and $ f(z)=\sum_{n=0}^{\infty}a_nz^n $. For $ N\in\mathbb{N} $, let $ t=\lfloor (N-1)/2 \rfloor $. Then 
	\begin{align}\label{e-11.88}
		|f(z)|+B_N(f,r)+sgn(t)\sum_{n=1}^{t}|a_n|^2\frac{r^N}{1-r}+\left(\frac{1}{1+|a_0|}+\frac{r}{1-r}\right)\sum_{n=t+1}^{\infty}|a_n|^2r^{2n}\leq 1
	\end{align} 
for $ |z|=r\leq R_N $, where $ R_N $ is the positive root of the equation $2(1+r)r^N-(1-r)^2=0$. The radius $ R_N $ is best possible. Moreover,
\begin{align}\label{e-11.99}
	|f(z)|^2+B_N(f,r)+sgn(t)\sum_{n=1}^{t}|a_n|^2\frac{r^N}{1-r}r+\left(\frac{1}{1+|a_0|}+\frac{r}{1-r}\right)\sum_{n=t+1}^{\infty}|a_n|^2r^{2n}\leq 1
\end{align}
for $ |z|=r\leq R^{\prime}_N $, where $ R^{\prime}_N $ is positive root of the equation $(1+r)r^N-(1-r)^2=0$. The radius $ R^{\prime}_N $ is best possible.
\end{thm}
\begin{rem}
In particular if $ N=1 $, then it is easy to see that $ R_1=\sqrt{5}-2 $ and $ R^{\prime}_1=1/3 $.
\end{rem}
 For $ p=1, 2 $, for $ f\in\mathcal{B} $, we define the functionals
\begin{align*}
	\mathcal{B}_{p, f}(z, r):=|f(z)|^p+B_1(f,r)+A(f_0,r).
\end{align*}In the context of Theorems \ref{Th-2.1} and \ref{TH-2.3}, the following result is obtained in \cite{Liu-Liu-Ponnusamy-BSM-2021} showing that the two constants can be improved for any individual function in $ \mathcal{B}. $
\begin{thm}\cite{Liu-Liu-Ponnusamy-BSM-2021}\label{th-2.5}
	Suppose that $ f\in\mathcal{B} $ and $ f(z)=\sum_{n=0}^{\infty}a_nz^n $. Then $ \mathcal{B}_{1, f}(z, r)\leq 1 $ for $ |z|=r\leq r_{a_0}=2/(3+|a_0|+\sqrt{5}(1+|a_0|))$. The radius $ r_{a_0} $ is best possible and $ r_{a_0}>\sqrt{5}-2 $. Moreover, $ \mathcal{B}_{2, f}(z, r)\leq 1 $ for $ |z|=r\leq r^{\prime}_{a_0} $, where $ r^{\prime}_{a_0} $ is the unique positive root of the equation
	\begin{align*}
		\left(1-|a_0|^3\right)r^3-(1+2|a_0|)r^2-2r+1=0.
	\end{align*}
 The radius $ r^{\prime}_{a_0} $ is best possible. Further, $ 1/3<r^{\prime}_{a_0}<1/(2+|a_0|) $.
\end{thm}
For recent developments on the Bohr-Rogosinski inequalities, we refer to the articles \cite{Aha-CMFT-2022,Allu-Arora-JMAA-2022,Allu-Halder-Pal-BDS-2022,Das-JMAA-2022} and references therein. However, we see that the quantities $ 1/(1+|a_0|)+r/(1-r) $ and $ \sum_{n=1}^{\infty}|a_n|^2r^{2n} $ for analytic functions in $ \mathcal{B} $ are analogous to $ 1/(1+|a_0|+|b_0|)+r/(1-r)=1+r/(1-r) $ (as $ a_0=0=b_0 $) and $ \sum_{n=2}^{\infty}(|a_n|+|b_n|)^2r^{2n} $, respectively, for harmonic functions given in \eqref{e-1.3}. The observations leads us to establish several harmonic analogues of the refined Bohr inequalities for certain class of harmonic mappings.
\subsection{Improved Bohr inequalities for the class $\mathcal{B}$}
Let $ f $ be holomorphic in $ \mathbb{D} $, and for $ 0<r<1 $, let $ \mathbb{D}_r=\{z\in\mathbb{C} : |z|<r\} $. The quantity $ S_r:=S_r(f) $ denoted as the planar integral has the following integral representation
\begin{align*}
	S_r:=\int_{\mathbb{D}_r}|f^{\prime}(z)|^2dA(z).
\end{align*}
Note that if $ f(z)=\sum_{n=0}^{\infty}a_nz^n $, then $ S_r:=\pi\sum_{n=1}^{\infty}n|a_n|^2r^{2n} $ (see \cite{Duren-1983}). It is well-known that if $ f $ is a univalent function, then  $ S_r $ is the area of the image of the sub-disk $ \mathbb{D}_r:=\{z\in\mathbb{D} : |z|<r\} $ under the mapping $ f $ (see \cite{Duren-1983}). The quantity $ S_r $ is non-negative and has been used extensively to study the improved versions of Bohr inequality and Bohr radius for the class of analytic functions (see e.g. \cite{Aha-Aha-CMFT-2023,Ismagilov-2020-JMAA,Kayumov-CRACAD-2018,Kayumov-Khammatova-JMAA-2021}) as well as harmonic mappings (see e.g. \cite{Ahamed-AMP-2021,Ahamed-CVEE-2021}). In the following, we recall some of the interesting results in this regard.  Recently, Kayumov and Ponnusamy \cite{Kayumov-CRACAD-2018} have obtained the following improved version of Bohr's inequality in terms of $ S_r $.
\begin{thm}\cite{Kayumov-CRACAD-2018}\label{th-1.12}
	Suppose that  $ f\in\mathcal{B} $ with $f(z)=\sum_{n=0}^{\infty} a_{n}z^{n}$. Then 
	\begin{equation}\label{e-11.13}
	B_0(f,r)+ \frac{16}{9} \left(\frac{S_{r}}{\pi}\right) \leq 1 \quad \mbox{for} \quad r \leq \frac{1}{3}
	\end{equation}
	and the numbers $1/3$, $16/9$ cannot be improved. Moreover, 
	\begin{equation}\label{e-11.14}
		|a_{0}|^{2}+B_1(f,r)+ \frac{9}{8} \left(\frac{S_{r}}{\pi}\right) \leq 1 \quad \mbox{for} \quad r \leq \frac{1}{2}
	\end{equation}
	and the numbers $1/2$, $9/8$ cannot be improved.
\end{thm}
Moreover, Ismagilov \emph{et al.} \cite{Ismagilov-2020-JMAA} investigated on the inequalities \eqref{e-11.13} and \eqref{e-11.14} of Theorem \ref{th-1.12} and obtained the following improved results.
\begin{thm}\cite{Ismagilov-2020-JMAA}\label{th-1.17}
	Suppose that  $ f\in\mathcal{B} $ with $f(z)=\sum_{n=0}^{\infty} a_{n}z^{n}$. Then 
	\begin{equation}\label{e-11.16}
	B_0(f,r)+\frac{16}{9}\left(\frac{S_r}{\pi}\right)+\lambda\left(\frac{S_r}{\pi}\right)^2\leq 1\;\; \mbox{for}\;\; r\leq\frac{1}{3},
	\end{equation}
	where
	\begin{equation*}
		\lambda=\frac{4(486-261a-324a^2+2a^3+30a^4+3a^5)}{81(1+a)^3(3-5a)}=18.6095...
	\end{equation*}
	and $ a\approx 0.567284 $ is the unique root in $ (0,1) $ of the equation  
	\begin{equation*}
		\phi(t)=-405+473r+402r^2+38r^3+3r^4+r^5.
	\end{equation*}
	The equality is achieved for the function $ f_a(z):=(a-z)/(1-az) $.
\end{thm}
\begin{thm}\cite{Ismagilov-2020-JMAA}\label{th-1.19}
	Suppose that  $ f\in\mathcal{B} $ with $f(z)=\sum_{n=0}^{\infty} a_{n}z^{n}$. Then 
	\begin{equation}\label{e-11.20}
		|f(z)|^2+B_1(f,r)+\frac{16}{9}\left(\frac{S_r}{\pi}\right)+\lambda\left(\frac{S_r}{\pi}\right)^2\leq 1\;\; \mbox{for}\;\; r\leq\frac{1}{3}
	\end{equation}
	where
	\begin{equation*}
		\lambda=\frac{-81+1044a+54a^2-116a^3-5a^4}{162(a+1)^2(2a-1)}=16.4618...
	\end{equation*}
	and $ a\approx 0.537869 $ is the unique root in $ (0,1) $ of the equation 
	\begin{equation*}
		-513+910t+80t^2+2t^3+t^4=0.
	\end{equation*}
	The equality is achieved for the function $ f_a $.
\end{thm}
We also recall here one result proved by Ismagilov \emph{et al.} \cite{Ismagilov-2020-JMAA} which is an improved version of the classical Bohr inequality, where $|a_0|$ is replaced by $|f(z)|.$
\begin{thm}	\cite{Ismagilov-2020-JMAA}\label{th-1.21}
	Suppose that  $ f\in\mathcal{B} $ with $f(z)=\sum_{n=0}^{\infty} a_{n}z^{n}$. Then 
	\begin{align}\label{e-11.22}
	|f(z)|+B_1(f,r)+p\left(\frac{S_r}{\pi}\right)\leq 1\;\; \mbox{for}\; |z|=r\leq \sqrt{5}-2,
	\end{align}
	where the constants $ r_0=\sqrt{5}-2\approx 0.236068 $ and $ p=2(\sqrt{5}-1) $ are sharp.
\end{thm}

Since, the quantity $ S_r/\pi $ is instrumental in the study of improved Bohr inequalities for its use to achieve the sharpness,   it is observed in \cite{Ismagilov-2021-JMS} the following inequality
\begin{align}\label{e-11.18}
	\dfrac{S_r}{\pi-S_r}\leq \dfrac{r^2(1-|a_0|^2)^2}{(1-r^2)(1-r^2|a_0|^4)}.
\end{align}
In view of the bound of the quantity $ S_r/(\pi-S_r) $, Ismagilov \emph{et al.} \cite{Ismagilov-2021-JMS} have investigated on Theorem \ref{th-1.12} and obtained the following sharp result.
\begin{thm} \label{th-1.18} \cite{Ismagilov-2021-JMS}
	Suppose that $ f(z)=\sum_{n=0}^{\infty}a_nz^n\in\mathcal{B} $. Then 
	\begin{align*}
		B_0(f,r)+ \frac{16}{9} \left(\frac{S_{r}}{\pi-S_r}\right) \leq 1 \quad \mbox{for} \quad r \leq \frac{1}{3},
	\end{align*}
	and the number $16/9$ cannot be improved. Moreover, 
	\begin{align*}
		|a_{0}|^{2}+B_1(f,r)+ \frac{9}{8} \left(\frac{S_{r}}{\pi-S_r}\right) \leq 1 \quad \mbox{for} \quad r \leq \frac{1}{2},
	\end{align*}
	and the number $9/8$ cannot be improved.
\end{thm}
 In recent developments, an refined version of the inequality \eqref{e-11.13} mentioned in \cite{Kayumov-CRACAD-2018} has been put forth. This improvement involves the substitution of the coefficient $|a_0|$ with $|f(z)|$, while employing the context of $S_r/(\pi-S_r)$ instead of $S_r/\pi$, as detailed below.
 \begin{thm}	\cite{Ismagilov-2020-JMAA}\label{th-2.11}
 	Suppose that  $ f\in\mathcal{B} $ with $f(z)=\sum_{n=0}^{\infty} a_{n}z^{n}$. Then 
 	\begin{align}\label{Eq-2.10}
 		|f(z)|+B_1(f,r)+2(\sqrt{5}-1)\left(\frac{S_r}{\pi-S_r}\right)\leq 1\;\; \mbox{for}\; |z|=r\leq \sqrt{5}-2,
 	\end{align}
 	where the constants $ r_0=\sqrt{5}-2\approx 0.236068 $ and $ 2(\sqrt{5}-1) $ are sharp.
 \end{thm}
 The motivation of the study Bohr phenomenon via proper combination for harmonic mappings is based on the discussions above and observation of the following remark concerning improved Bohr inequality.
 \begin{rem}\label{Rem-2.1}
 	Ismagilov \emph{et al.} remarked in Theorem \ref{th-1.17} that for any function $F : [0, \infty)\to [0, \infty)$ such that $F(t)>0$ for $t>0$, there exists analytic function $f : \mathbb{D}\to\mathbb{D}$ for which the inequality
 	\begin{align*}
 		\sum_{n=0}^{\infty}|a_n|r^n+\frac{16}{9}\left(\frac{S_r}{\pi}\right)+\lambda\left(\frac{S_r}{\pi}\right)^2+F(S_r)\leq 1\;\; \mbox{for}\;\; r\leq\frac{1}{3}
 	\end{align*}
 	is false, there $\lambda$ is given in Theorem \ref{th-1.17}.
 \end{rem}
 Considering Remark \ref{Rem-2.1}, it is noteworthy to observe that augmenting a non-negative quantity with the Bohr inequality does not yield the desired inequality for the class $\mathcal{B}$. This observation motivates us to pose the following question for further study on the topic.
 \begin{ques}\label{Q1}
 	What can be deduced about the refine harmonic analogue of Theorem \ref{TH-2.3} by incorporating a general power of $ |f(z)|^p $ for $ p\in\mathbb{N} $ and adding  higher power terms of $S_r/\pi$ for the class $\mathcal{P}^{0}_{\mathcal{H}}(M)?$ Consequently, can we obtain the harmonic analogue of the Theorems \ref{th-2.5} to \ref{th-1.21} for the class $\mathcal{P}^{0}_{\mathcal{H}}(M)?$
 \end{ques}
In the past few years, there has been a noteworthy investigation into substituting $S_r/(\pi-S_r)$ for $S_r/\pi$ in mazorent series and exploring the determination of the Bohr radius in this modified framework, making it a compelling subject in geometric function theory. Motivated from the work of Ismagilov \emph{et al.} \cite{Ismagilov-2021-JMS}, it is natural to raise the following question also.
\begin{ques}\label{Q2}
 What can be deduced about the refine harmonic counterpart of Theorem \ref{TH-2.3} by incorporating a general power of $ |f(z)|^p $ for $ p\in\mathbb{N} $ and adding higher power terms of $S_r/(\pi-S_r)$ instead of the term $S_r/\pi$ for the class $\mathcal{P}^{0}_{\mathcal{H}}(M)?$ Moreover, can we obtain the harmonic analogue of the Theorems \ref{th-1.18} and \ref{th-2.11} for the class $\mathcal{P}^{0}_{\mathcal{H}}(M)?$
\end{ques}
The organization of this paper is as follows. In Section 3, we state the main results to address Question \ref{Q1} and Question \ref{Q2} for the close to convex class  $\mathcal{P}^{0}_{\mathcal{H}}(M)$ of harmonic mappings. Furthermore, we shall discuss some consequences of the main results showing Bohr radii by graphs and tables for the class $ \mathcal{P}^{0}_{\mathcal{H}}(M) $ for different values of $M.$ The section 4 contains the proof of our main results.
 
\section{Main results and corollaries}
Before stating the main results, we recall here the definition of dilogarithm $Li_2(z)$ which will play a crucial role in our study. In fact, the polylogarithm $Li_n(z)$ is  defined as
\begin{align*}
	Li_n(z)=\sum_{k=1}^{\infty}\frac{z^k}{k^n}=z+\frac{z^2}{2^n}+\frac{z^3}{3^n}+\cdots
\end{align*}
The polylogarithm is a speical function $Li_n(z)$ of order $n$ and the argument $z$. 
 For, $n=2,$  we get the dilogarithm $Li_2(z)$ defined by
\begin{equation*}
	{\rm Li}_2(z)=\sum_{k=1}^{\infty}\frac{z^k}{k^2}\;\; \mbox{for}\;\; |z|<1.
\end{equation*}
The following properties of dilogarithm hold
\begin{equation*}
	{\rm Li}_2(r)+{\rm Li}_2(1-r)=\frac{\pi^2}{6}-\log r\log (1-r);
\end{equation*}
\begin{align*}
	{\rm Li}_2(1-r^2)\rightarrow 0\;\; \mbox{and}\;\; \frac{\log(1/r)}{1-r}\rightarrow 1\;\; \mbox{as}\;\; r\rightarrow 1.
\end{align*} 
In particular, the analytic continuation of the dilogarithm is given by
\begin{equation*}
	{\rm Li}_2(z)=-\int_{0}^{z}\log (1-u)\frac{du}{u}\;\; \mbox{for}\;\; z\in\mathbb{C}\setminus [1,\infty).
\end{equation*}

For our further discussions, we need to introduce the following notations also.
Let $P_q(w)$ be a polynomial defined as
\begin{align*}
	P_q(w):=\lambda_1w+\lambda_2w^2+\cdots+\lambda_qw^q\; \mbox{for}\; q\in\mathbb{N}\cup\{0\},
\end{align*}where $\lambda_j\in\mathbb{R}_{\geq 0}:=\{x\in\mathbb{R} : x\geq 0\} $ for all  $j=0,1,\cdots q.$\vspace{2mm}

Since our aim is to study Bohr Phenomenon both refined and improved inequalities, we consider $	H^{f}_{\beta, \mu, \lambda, m, N}(r)$ as suitable combination for harmonic mappings belong to the class $\mathcal{P}^0_H(M).$ Moreover, we define the functionals $	H^{f}_{\beta, \mu, \lambda, m, N}(r)$ and $G^{f,\lambda_1,\lambda_2,\cdots\lambda_q}_{\beta, \mu, \lambda, m, N}(r)$ are as follows:
\begin{align*}
	H^{f}_{\beta, \mu, \lambda, m, N}(r):&=\beta|f(z)|^m+\sum_{n=N}^{\infty}\left(|a_n|+|b_n|\right)r^n+\mu\; sgn(t)\sum_{n=1}^{t}\left(|a_n|+|b_n|\right)^2\frac{r^N}{1-r}\\&\quad+\lambda\left(1+\frac{r}{1-r}\right)\sum_{n=t+1}^{\infty}\left(|a_n|+|b_n|\right)^2r^{2n} 
\end{align*}and
\begin{align*}
		G^{f,\lambda_1,\lambda_2,\cdots\lambda_q}_{\beta, \mu, \lambda, m, N}(r):&=H^{f}_{\beta, \mu, \lambda, m, N}(r)+P_q\bigg(\frac{S_r}{\pi}\bigg),
\end{align*}
where $ N $ be a positive integer and $\beta\in\mathbb{R}_{\geq 0}$. 
\begin{rem}
One fact worth mentioning here is that in the proof of main results, a term $ \sum_{n=2}^{t}\frac{r^{2n}}{n^2(n-1)^2} $ will appear in our computation, hence following range of the summation (i.e. $ t\geq 2 $), if $ N=1, 2 $, then $ t=\lfloor (N-1)/2 \rfloor=0 $ and hence $ sgn(t)=0. $ Further, if $ N=3, 4 $, then we see that $ t=\lfloor (N-1)/2 \rfloor=1 $ and if $ N\geq 5, $ then $ t\geq 2 $ with $ sgn(t)=1 $. In view of this observations, to serve our purpose, in this paper, in all the main results, we will consider $ N\geq 5 $. The possible situation for the cases $ N=1, 2, 3, 4 $, we give corollary of the corresponding main result.
\end{rem}
 
The main aim of this paper is to establish several Bohr inequalities, finding the corresponding sharp radius for the class $ \mathcal{P}^{0}_{\mathcal{H}}(M) $, many other properties of the class are established by Ghosh and Vasudevarao(See \cite{Ghosh-Vasudevarao-BAMS-2020})
$$\mathcal{P}^{0}_{\mathcal{H}}(M)=\{f=h+\overline{g} \in \mathcal{H}_{0}: \real (zh^{\prime\prime}(z))> -M+|zg^{\prime\prime}(z)|, \; z \in \mathbb{D}\; \mbox{and }\; M>0\}.$$

However, the following lemmas (see \cite{Ghosh-Vasudevarao-BAMS-2020}) play key roles to prove the main results. It provides the coefficient bounds and the growth estimates for functions in the class  $ \mathcal{P}^{0}_{\mathcal{H}}(M)$.
\begin{lem} \label{lem-2.19}
	Let $f=h+\overline{g}\in \mathcal{P}^{0}_{\mathcal{H}}(M)$ be given by \eqref{e-1.3} for some $M>0$. Then for $n\geq 2,$ 
	(i) $\displaystyle |a_n| + |b_n|\leq \frac {2M}{n(n-1)}; $
	(ii) $\displaystyle ||a_n| - |b_n||\leq \frac {2M}{n(n-1)};$
	(iii) $\displaystyle |a_n|\leq \frac {2M}{n(n-1)}.$	The inequalities  are sharp with extremal function   $f_M$ given by 
	$f_M^{\prime}(z)=1-2M\, \ln\, (1-z) .$	
\end{lem}
\begin{lem}\label{lem-2.20}
	Let $f \in \mathcal{P}^{0}_{\mathcal{H}}(M)$ be given by \eqref{e-1.3}. Then 
	\begin{equation} \label{e-2.15}
		|z| +2M \sum\limits_{n=2}^{\infty} \dfrac{(-1)^{n-1}|z|^{n}}{n(n-1)} \leq |f(z)| \leq |z| + 2M \sum\limits_{n=2}^{\infty} \dfrac{|z|^{n}}{n(n-1)}.
	\end{equation}
	Both  inequalities are sharp for the function $f_{M}$ given by $f_{M}(z)=z+ 2M \sum\limits_{n=2}^{\infty} \dfrac{z^n}{n(n-1)}.
	$
\end{lem}
The following is our first main result, which addresses the Question \ref{Q1} and estimates the generalized refine Bohr-Rogosinski inequality for  the class  $ \mathcal{P}^{0}_{\mathcal{H}}(M). $
\begin{thm}\label{th-3.1}
Let $ f\in \mathcal{P}^{0}_{\mathcal{H}}(M) $ be given by \eqref{e-1.3} and $ 0\leq M <1/(2(\ln 4-1))$. Then for $ \mu,\; \lambda,\lambda_j\in\mathbb{R}_{\geq 0}:=\{x\in\mathbb{R} : x\geq 0\} $,where $j=1,2,\cdots q $ and $ N\geq 5 $, we have $ G^{f,\lambda_1,\lambda_2,\cdots\lambda_q}_{\beta, \mu, \lambda, m, N}(r)\leq {d}\left(f(0),\partial \mathbb{D}\right) $ for $ |z|=r\leq R^{\lambda_1,\lambda_2,\cdots\lambda_q}_{\beta,\mu, \lambda,m,N}(M)  $, where $ R^{\lambda_1,\lambda_2,\cdots\lambda_q}_{\beta,\mu, \lambda,m,N}(M) $ is the unique root  in $ (0,1) $ of the equation $ 	\Phi^{\lambda_1,\lambda_2,\cdots\lambda_q}_{\beta,\mu, \lambda,m,N}(r)=0, $ where
\begin{align}\label{eq-2.4}
\Phi^{\lambda_1,\lambda_2,\cdots\lambda_q}_{\beta,\mu, \lambda,m,N}(r):\nonumber&= \beta\left(G_{M}(r)\right)^m+2M\left(r+(1-r)\ln(1-r)-\sum_{n=2}^{N-1}\frac{r^n}{n(n-1)}\right)+\mu \Phi^N_{M,t}(r)\\&\quad+4M^2\lambda\left(1+\frac{r}{1-r}\right) G_{2,t}(r)+P_q(F_M(r))-1-2M\left(1-2\ln 2\right)
\end{align}
	and
	\[
	\begin{cases}
G_{M}(r):=\displaystyle r+2M\sum_{n=2}^{\infty}\frac{r^n}{n(n-1)}=r+2M\left(r+(1-r)\ln (1-r)\right);\vspace{1.4mm}\\
\Phi^N_{M,t}(r):=\displaystyle \frac{4M^2r^N}{1-r}sgn(t)\sum_{n=1}^{t}\frac{1}{n^2(n-1)^2};\vspace{1.4mm}\\
G_{t}(r):=\displaystyle\left(r^2+1\right){\rm Li_2\left(r^2\right)}+2\left(r^2-1\right)\ln\left(1-r^2\right)-3r^2-\sum_{n=2}^{t}\frac{r^{2n}}{n^2(n-1)^2};\vspace{1.4mm}\\
 F_{M}(r):=r^2+4M[r^2Li_2(r^2)-(r^2+(1-r^2)\log(1-r^2)).
	\end{cases}
	\]
	The constant $ R^{\lambda_1,\lambda_2,\cdots\lambda_q}_{\beta,\mu, \lambda,m,N}(M) $ is best possible.
\end{thm}
For a sake of simplification, we have used the following assumptions  in order to present certain consequences of Theorem \ref{th-3.1} here:
\[
\begin{cases}
\mathcal{J}^M_1(r):=r^2+4M^2[\left(1+r^2\right){\rm Li_2\left(r^2\right)}-2\left(1-r^2\right)\ln\left(1-r^2\right)-3r^2];\\
\mathcal{J}^M_2(m, r):=(G_M(r))^m-1-2M(1-2\ln 2);\\
\mathcal{J}_3(r):=r+(1-r)\ln(1-r).
\end{cases}
\]
In fact, the following corollary of Theorem \ref{th-3.1}, we obtain the following result in which the cases for $ N=1, 2, 3, 4 $ are discussed completely.
\begin{cor} Let $ f\in \mathcal{P}^{0}_{\mathcal{H}}(M) $ be given by \eqref{e-1.3} and $ 0\leq M <1/(2(\ln 4-1))$, $ \mu,\; \lambda\in\mathbb{R}_{\geq 0} $.
\begin{enumerate}
\item[(i)] If $ N=1 $, then $ G^{f,\lambda_1,\lambda_2,\cdots\lambda_q}_{\beta, \mu, \lambda, m, 1}(r)\leq {d}\left(f(0),\partial \mathbb{D}\right) $
for $ |z|=r\leq R^{\lambda_1,\lambda_2,\cdots\lambda_q}_{\beta,\mu, \lambda,m,1}(M) $, where $ R^{\lambda_1,\lambda_2,\cdots\lambda_q}_{2,\mu, \lambda,m,1}(M)  $ is the unique root in $ (0,1) $ of the equation 
\begin{align*}
r+2M\mathcal{J}_3(r)+\mathcal{J}^M_2(m, r)+\lambda\left(1+\frac{r}{1-r}\right)\mathcal{J}^M_1(r)+P_q(F_{M}(r))=0.
\end{align*}
\item[(ii)] If $ N=2 $, then $ G^{f,\lambda_1,\lambda_2,\cdots\lambda_q}_{\beta, \mu, \lambda, m, 2}(r)\leq {d}\left(f(0),\partial \mathbb{D}\right) $
for $ |z|=r\leq R^{\lambda_1,\lambda_2,\cdots\lambda_q}_{\beta,\mu, \lambda,m,2}(M) $, where $ R^{\lambda_1,\lambda_2,\cdots\lambda_q}_{\beta,\mu, \lambda,m,2}(M) $ is the unique root in $ (0,1) $ of the equation 
\begin{align*}
2M\mathcal{J}_3(r)+\mathcal{J}^M_2(m, r)+\lambda\left(1+\frac{r}{1-r}\right)\mathcal{J}^M_1(r)+P_q(F_{M}(r))=0.
\end{align*}
\item[(iii)] If $ N=3 $, then $ G^{f,\lambda_1,\lambda_2,\cdots\lambda_q}_{\beta, \mu, \lambda, m, 3}(r)\leq {d}\left(f(0),\partial \mathbb{D}\right) $
for $ |z|=r\leq R^{\lambda_1,\lambda_2,\cdots\lambda_q}_{\beta,\mu, \lambda,m,3}(M)$, where $ R^{\lambda_1,\lambda_2,\cdots\lambda_q}_{\beta,\mu, \lambda,m,3}(M) $ is the unique root in $ (0,1) $ of the equation 
\begin{align*}
&2M\left(\mathcal{J}_3(r)-\frac{r^2}{2}\right)+\mathcal{J}^M_2(m, r)+\mu\frac{r^3}{1-r}+\lambda\left(1+\frac{r}{1-r}\right)\left(\mathcal{J}^M_1(r)-r^2\right)\\&\quad+P_q(F_{M}(r))=0.
\end{align*}
\item[(iv)] If $ N=4 $, then $ G^{f,\lambda_1,\lambda_2,\cdots\lambda_q}_{\beta, \mu, \lambda, m, 4}(r)\leq {d}\left(f(0),\partial \mathbb{D}\right) $
for $ |z|=r\leq R^{\lambda_1,\lambda_2,\cdots\lambda_q}_{\beta,\mu, \lambda,m,4}(M) $, where $ R^{\lambda_1,\lambda_2,\cdots\lambda_q}_{\beta,\mu, \lambda,m,4}(M)  $ is the unique root in $ (0,1) $ of the equation 
\begin{align*}
&2M\left(\mathcal{J}_3(r)-\frac{r^2}{2}-\frac{r^3}{6}\right)+\mathcal{J}^M_2(m, r)+\mu\frac{r^4}{1-r}+\lambda\left(1+\frac{r}{1-r}\right)\left(\mathcal{J}^M_1(r)-r^2\right)\quad\\&\quad+P_q(F_{M}(r))=0.
\end{align*}
\end{enumerate}
The constants $ R^{\lambda_1,\lambda_2,\cdots\lambda_q}_{\beta,\mu, \lambda,m,1}(M)$, $R^{\lambda_1,\lambda_2,\cdots\lambda_q}_{\beta,\mu, \lambda,m,2}(M) $, $ R^{\lambda_1,\lambda_2,\cdots\lambda_q}_{\beta,\mu, \lambda,m,3}(M) $ and $ R^{\lambda_1,\lambda_2,\cdots\lambda_q}_{\beta,\mu, \lambda,m,4}(M) $ all are best possible.
\end{cor}
For $m=1,2$, we define a functional $G_{m}^f$ 
 \begin{align*}
 	G_{m}^{f}:=|f(z)|^m+r+\sum_{n=2}^{\infty}\left(|a_n|+|b_n|\right)r^n+\left(1+\frac{r}{1-r}\right)\sum_{n=1}^{\infty}\left(|a_n|+|b_n|\right)^2r^{2n}.
 \end{align*}
The following corollary is a harmonic analogue of the Theorem \ref{th-2.5} for the class $\mathcal{P}^{0}_{\mathcal{H}}(M).$ A part of Question \ref{Q1} is thus answered.\vspace{2mm}
\begin{cor}	Let $ f\in \mathcal{P}^{0}_{\mathcal{H}}(M) $ be given by \eqref{e-1.3} and $ 0\leq M <1/(2(\ln 4-1))$. Then 
	\begin{enumerate}
		\item [(i)] $G_{1}^f\leq {d}\left(f(0),\partial \mathbb{D}\right)$ for $|z|=r\leq R^{0,0,\cdots0}_{1,\mu, 1,1,1}(M) $, where $R^{0,0,\cdots0}_{1,\mu, 1,1,1}(M) $ is the unique root in $ (0,1) $ of the equation 
		\begin{align}\label{eq-3.8}
			2r-1+2M\big(2r-1+\ln4+(1-r)\ln(1-r)^2\big)+\bigg(1+\frac{r}{1-r}\bigg)\mathcal{J}^M_1(r)=0
		\end{align}where the constant $R^{0,0,\cdots0}_{1,\mu, 1,1,1}(M)$ is best possible.
		
		\item[(ii)] $G_{2}^f\leq {d}\left(f(0),\partial \mathbb{D}\right)$  for $|z|=r\leq R^{0,0,\cdots0}_{1,\mu, 1,2,1}(M) $, where $R^{0,0,\cdots0}_{1,\mu, 1,2,1}(M) $ is the unique root in $ (0,1) $ of the equation 
		\begin{align}\label{eq-3.9}
			\left(G_M(r)\right)^2+r-1+2M\big((1-r)\big(\ln(1-r)-1\big)+\ln4+\big)+\bigg(1+\frac{r}{1-r}\bigg)\mathcal{J}^M_1(r)=0
		\end{align}where the constant $R^{0,0,\cdots0}_{1,\mu, 1,2,1}(M)$ is best possible.
	\end{enumerate} 
		\begin{figure}
		\begin{center}
			\includegraphics[width=0.5
			\linewidth]{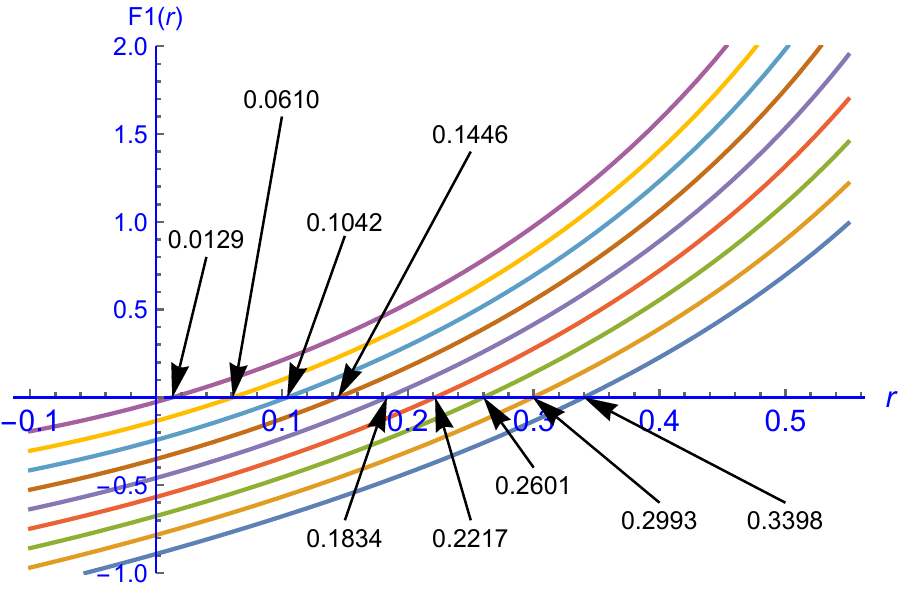}\;\;\;\;\;\includegraphics[width=0.5\linewidth]{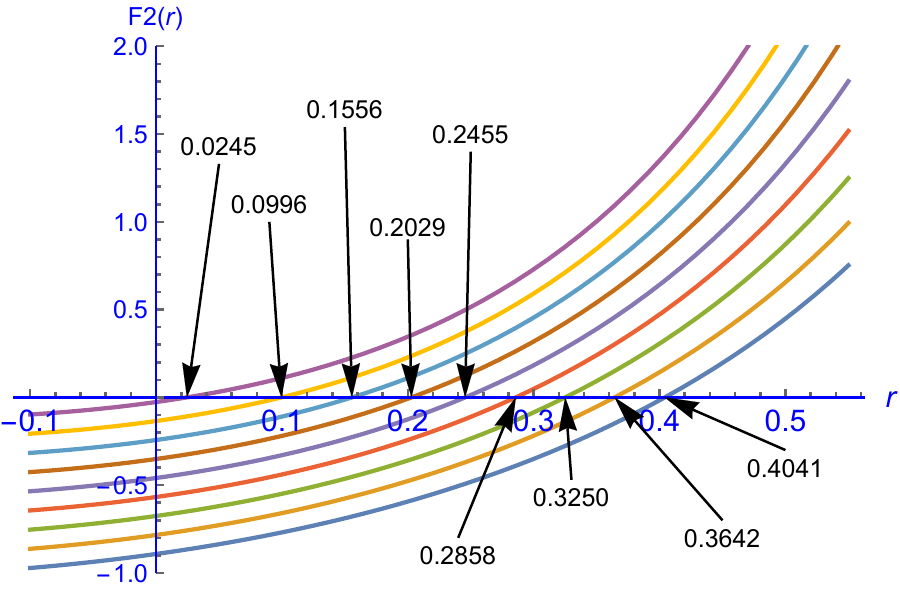}
		\end{center}
		\caption{The figure exhibits the roots $R^{0,0,\cdots0}_{1,\mu, 1,1,1}(M)$ of equation \eqref{eq-3.8} and the roots  $R^{0,0,\cdots0}_{1,\mu, 1,2,1}(M)$  of   \eqref{eq-3.9} respectively for different values of $M$ as defined in Table \ref{t-1}}
	\end{figure}
	\begin{table}[ht]
		\centering
		\begin{tabular}{|l|l|l|l|l|l|l|l|l|l|}
			\hline
			$M$& $0.14 $&$0.28 $& $0.42 $& $0.56$& $0.70 $&$0.84 $&$ 0.98 $ &$1.12$& $1.26$  \\
			\hline
			$ R^{0,0,\cdots0}_{1,\mu, 1,1,1}(M)$& $0.3398$&$0.2993 $& $0.2601$& $0.2217$& $0.1834 $&$0.1446 $& $0.1042$& $0.0610$ &$ 0.0129 $\\
			\hline
			$R^{0,0,\cdots0}_{1,\mu, 1,2,1}(M)$& $0.4041$&$0.3642 $& $0.3250$& $0.2858$& $0.2455 $&$0.2029 $& $0.1556$& $0.0996$ &$ 0.0245 $\\
			\hline
		\end{tabular}\vspace{2.5mm}
		\caption{The table exhibits the roots of the equations \eqref{eq-3.8} and \eqref{eq-3.9} for different values of $M.$}
		\label{t-1}
	\end{table}
\end{cor}
We obtain the next corollary as a  harmonic analogue of the Theorem \ref{th-1.12} for the class $\mathcal{P}^{0}_{\mathcal{H}}(M).$ Thus a part of Question \ref{Q1} is answered.
\begin{cor}	Let $ f\in \mathcal{P}^{0}_{\mathcal{H}}(M) $ be given by \eqref{e-1.3} and $ 0\leq M <1/(2(\ln 4-1))$. Then  we have
	\begin{enumerate}
		\item[(a)] \begin{align*}
			r+\sum_{n=2}^{\infty}\left(|a_n|+|b_n|\right)r^n+\frac{16}{9}\bigg(\frac{S_r}{\pi}\bigg)\leq {d}\left(f(0),\partial \mathbb{D}\right)
		\end{align*}
		for $ |z|=r\leq R^{\frac{16}{9},0,\cdots0}_{0,\mu, 0,m,1}(M) $, where $R^{\frac{16}{9},0,\cdots0}_{0,\mu, 0,m,1}(M)  $ is the unique root in $ (0,1) $ of the equation 
		\begin{align}\label{EQN-3.5}
			r-1+2M\big((1-r)\big(\log(1-r)-1\big)+\ln4\big)+\dfrac{16}{9}\big(F_{M}(r)\big)=0.
		\end{align}The constant $R^{\frac{16}{9},0,\cdots0}_{0,\mu, 0,1,1}(M)$ is best possible. Moreover
		\item[(b)]
		\begin{align*}
			r+\sum_{n=2}^{\infty}\left(|a_n|+|b_n|\right)r^n+\frac{9}{8}\bigg(\frac{S_r}{\pi}\bigg)\leq {d}\left(f(0),\partial \mathbb{D}\right)
		\end{align*}
		for $ |z|=r\leq R^{\frac{9}{8},0,\cdots0}_{0,\mu, 0,m,1}(M) $, where $R^{\frac{9}{8},0,\cdots0}_{0,\mu, 0,m,1}(M) $ is the unique root in $ (0,1) $ of the equation
		\begin{align}\label{EQN-3.6}
			r-1+2M\big((1-r)\big(\log(1-r)-1\big)+\ln4\big)+\dfrac{9}{8}\big(F_{M}(r)\big)=0.
		\end{align}
		The constant $R^{\frac{9}{8},0,\cdots0}_{0,\mu, 0,2,1}(M)$ is best possible.
	\end{enumerate}
\end{cor}
\begin{figure}
	\begin{center}
		\includegraphics[width=0.5
		\linewidth]{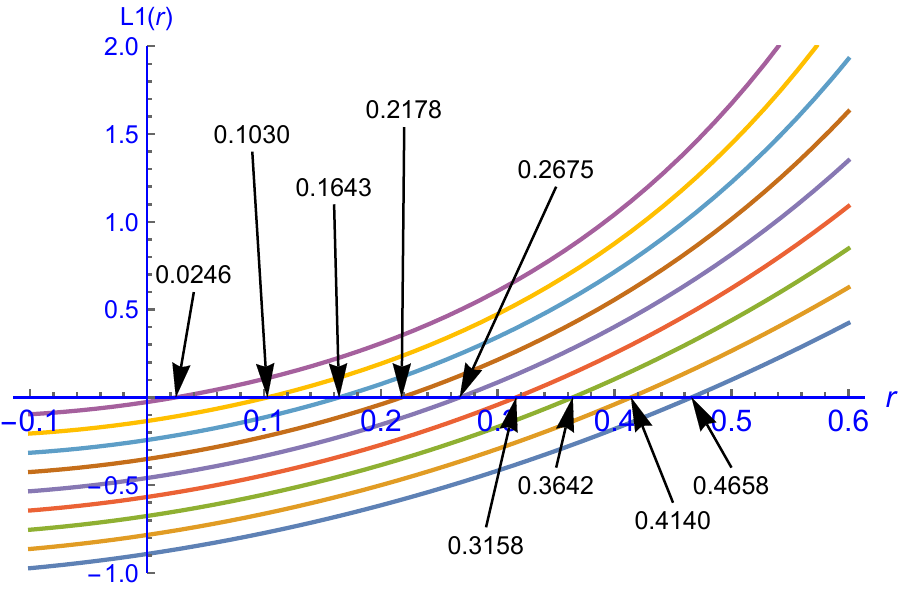 }\;\;\;\;\;\includegraphics[width=0.5\linewidth]{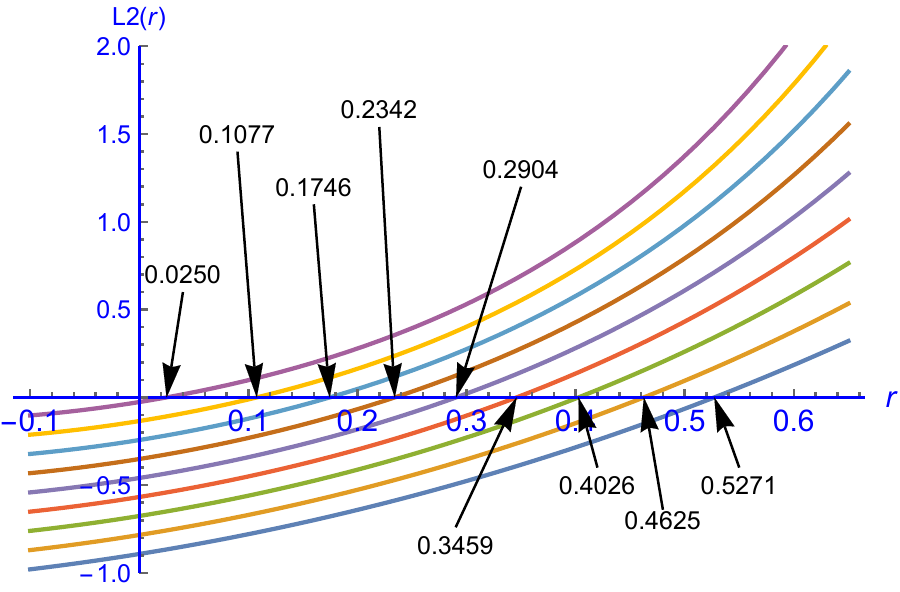}
	\end{center}
	\caption{The figure exhibits the roots $R^{\frac{16}{9},0,\cdots0}_{0,\mu, 0,m,1}(M)$ of equation \eqref{EQN-3.5} and the roots  $R^{\frac{9}{8},0,\cdots0}_{0,\mu, 0,m,1}(M)$  of   \eqref{EQN-3.6} respectively for different values of $M$ as defined in Table \ref{T-1}}
\end{figure}
\begin{table}[ht]
	\centering
	\begin{tabular}{|l|l|l|l|l|l|l|l|l|l|}
		\hline
		$M$& $0.14 $&$0.28 $& $0.42 $& $0.56$& $0.70 $&$0.84 $&$ 0.98 $ &$1.12$& $1.26$  \\
		\hline
		$R^{\frac{16}{9},0,\cdots0}_{0,\mu, 0,m,1}(M)$& $0.4658$&$0.4140 $& $0.3642$& $0.3158$& $0.2675 $&$0.2178$& $0.1643$& $0.1030$ &$ 0.0246 $\\
		\hline
		$ R^{\frac{9}{8},0,\cdots0}_{0,\mu, 0,m,1}(M)$& $0.5271$&$0.4625 $& $0.4026$& $0.3459$& $0.2904 $&$0.2342 $& $0.1746$& $0.1077$ &$ 0.0250 $\\
		\hline
	\end{tabular}\vspace{2.5mm}
	\caption{The table exhibits the roots of the equations \eqref{EQN-3.5} and \eqref{EQN-3.6} for different values of $M.$}
	\label{T-1}
\end{table}

  The following corollary is a harmonic analogue of the Theorem \ref{th-1.21} for the class $\mathcal{P}^{0}_{\mathcal{H}}(M).$ A part of Question \ref{Q1} is thus answered.
\begin{cor}Let $ f\in \mathcal{P}^{0}_{\mathcal{H}}(M) $ be given by \eqref{e-1.3} and $ 0\leq M <1/(2(\ln 4-1))$. Then we have
\begin{align*}
|f(z)|+r+\sum_{n=2}^{\infty}\left(|a_n|+|b_n|\right)r^n+2(\sqrt{5}-1)\bigg(\frac{S_r}{\pi}\bigg)\leq {d}\left(f(0),\partial \mathbb{D}\right)
\end{align*}  
for $ |z|=r\leq R^{2(\sqrt{5}-1),0,\cdots0}_{1, \mu, 0, 1, 1}(M) $, where $R^{2(\sqrt{5}-1),0,\cdots0}_{1, \mu, 0, 1, 1}(M)$ is the unique root  in $ (0,1) $ of the equation 
\begin{align}\label{eq-3.11}
r-1+2M\big(r-1+\ln4+(1-r)\ln(1-r)\big)+2(\sqrt{5}-1)(F_{M}(r))=0.
\end{align} 
The constant $R^{2(\sqrt{5}-1),0,\cdots0}_{1, \mu, 0, 1, 1}(M)$ is best possible.
\begin{figure}
\begin{center}
    \includegraphics[width=0.5\linewidth]{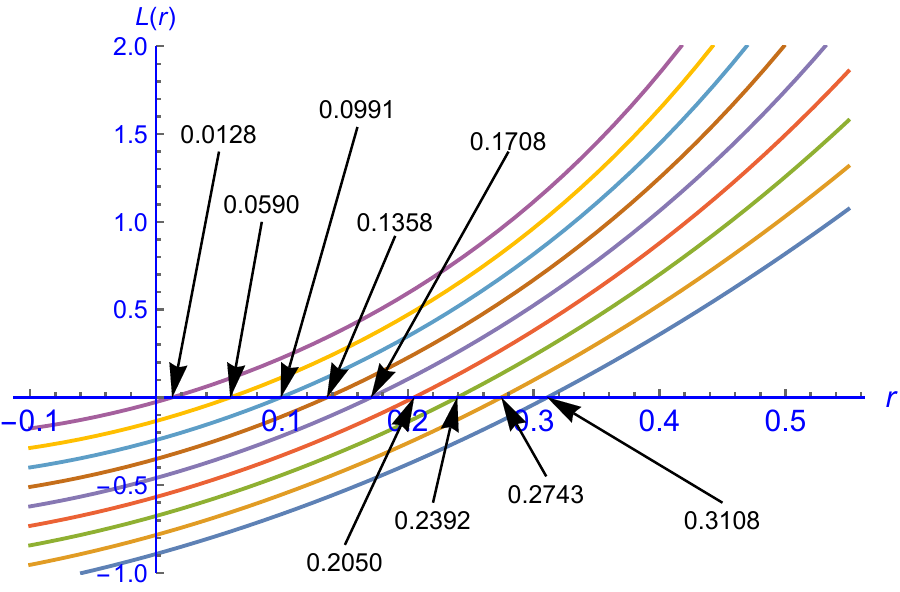}\;\;\;\;\includegraphics[width=0.5\linewidth]{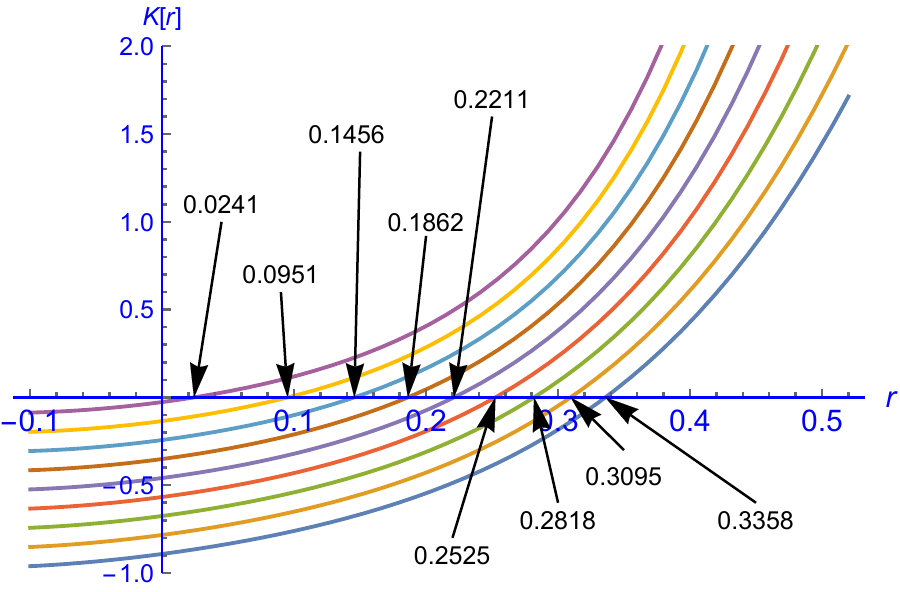}
		\end{center}
		\caption{The figure exhibits the roots of equation \eqref{eq-3.11} and roots of the equation \eqref{Eq-3.12}  as defined in Table \ref{t-4} and Table \ref{t-5} respectively.}
	\end{figure}
	\begin{table}[ht]
		\centering
		\begin{tabular}{|l|l|l|l|l|l|l|l|l|l|}
			\hline
			$M$& $0.14 $&$0.28 $& $0.42 $& $0.56$& $0.70 $&$0.84 $&$ 0.98 $ &$1.12$& $1.26$  \\
			\hline
			$R^{2(\sqrt{5}-1),0,\cdots0}_{1, \mu, 0, 1, 1}(M)$& $0.3108$&$0.2743 $& $0.2392$& $0.2050$& $0.1708 $&$0.1358 $& $0.0991$& $0.0590$ &$ 0.0128 $\\
			\hline
		\end{tabular}\vspace{2.5mm}
		\caption{The table exhibits the roots of the equation \eqref{eq-3.11}  for different values of $M.$}
		\label{t-4}
	\end{table}
\end{cor}
The following corollary is a harmonic analogue of the Theorem \ref{th-1.19} for the class $\mathcal{P}^{0}_{\mathcal{H}}(M).$ Consequently, a part of Question \ref{Q1} is answered.
\begin{cor}Let $ f\in \mathcal{P}^{0}_{\mathcal{H}}(M) $ be given by \eqref{e-1.3} and $ 0\leq M <1/(2(\ln 4-1))$. Then we have
	\begin{align*}
		|f(z)|^2+r+\sum_{n=2}^{\infty}\left(|a_n|+|b_n|\right)r^n+\frac{16}{9}\bigg(\frac{S_r}{\pi}\bigg)+\lambda_2\bigg(\dfrac{S_r}{\pi}\bigg)^2\leq {d}\left(f(0),\partial \mathbb{D}\right) 
	\end{align*}  
	for $ |z|=r\leq R^{16/9,\lambda_2,0,\cdots0}_{1, \mu, 0, 2, 1}(M)  $, where $R^{16/9,\lambda_2,0,\cdots0}_{1, \mu, 0, 2, 1}(M) $ is the unique root in $ (0,1) $ of the equation
	\begin{align}\label{Eq-3.12}
		\left(G_M(r)\right)^2+r-1+2M\big(r-1+\ln4+(1-r)\ln(1-r)\big)+\dfrac{16}{9}(F_{M}(r))+\lambda_2(F_{M}(r))^2=0,
	\end{align}
	where $\lambda_2=\lambda$ is as in Theorem \ref{th-1.19}.	The constant $R^{16/9,\lambda_2,0,\cdots0}_{1, \mu, 0, 2, 1}(M)$ is best possible.
	\begin{table}[ht]
		\centering
		\begin{tabular}{|l|l|l|l|l|l|l|l|l|l|}
			\hline
			$M$& $0.14 $&$0.28 $& $0.42 $& $0.56$& $0.70 $&$0.84 $&$ 0.98 $ &$1.12$& $1.26$  \\
			\hline
			$R^{16/9,\lambda_2,0,\cdots0}_{1, \mu, 0, 2, 1}(M)$& $0.3358$&$0.3095 $& $0.2818$& $0.2525$& $0.2211 $&$0.1862 $& $0.1456$& $0.0951$ &$ 0.0241 $\\
			\hline
		\end{tabular}\vspace{2.5mm}
		\caption{The table exhibits the roots of the equation \eqref{Eq-3.12}  for different values of $M$, where $\lambda_2=\lambda$ is as in Theorem \ref{th-1.19}.}
		\label{t-5}
	\end{table}
\end{cor}
Inspired by the idea of considering the quantity in Theorem \ref{th-1.18}, our next aim is to explore sharp Bohr inequalities for the class $\mathcal{P}^{0}_{\mathcal{H}}(M)$. Henceforth, we define the following functional
\begin{align}\label{eq-3.7}
	\mathcal{S}^{f,\lambda_1,\lambda_2,\cdots\lambda_q}_{\beta, \mu, \lambda, m, N}(r):&=H^{f}_{\beta, \mu, \lambda, m, N}(r)+P_q\bigg(\frac{S_r}{\pi-S_r}\bigg).
\end{align}
To address the Question \ref{Q2} and serve our purpose, we obtain the following result for the class $\mathcal{P}^{0}_{\mathcal{H}}(M)$.
\begin{thm}\label{th-3.2}
	Let $ f\in \mathcal{P}^{0}_{\mathcal{H}}(M) $ be given by \eqref{e-1.3} and $ 0\leq M <1/(2(\ln 4-1))$. Then for $ \mu,\; \lambda,\lambda_j\in\mathbb{R}_{\geq 0}:=\{x\in\mathbb{R} : x\geq 0\} $,where $j=1,2,\cdots q $ and $ N\geq 5 $, we have 
	\begin{align*}
		\mathcal{S}^{f,\lambda_1,\lambda_2,\cdots\lambda_q}_{\beta,\mu, \lambda, m, N}(r)\leq {d}\left(f(0),\partial \mathbb{D}\right) \;\mbox{for}\;   |z|=r\leq R^{*, \lambda_1,\lambda_2,\cdots\lambda_q}_{\beta,\mu, \lambda,m,N}(M)  
	\end{align*}
	 where $ R^{*, \lambda_1,\lambda_2,\cdots\lambda_q}_{\beta,\mu, \lambda,m,N}(M) $ is the unique root of the equation $ 	\Phi^{*, \lambda_1,\lambda_2,\cdots\lambda_q}_{\beta,\mu, \lambda,m,N}(r)=0, $  where
	\begin{align*}
		\Phi^{*,\lambda_1,\lambda_2,\cdots\lambda_q}_{\beta,\mu, \lambda,m,N}(r):\nonumber&= \beta\left(G_{M}(r)\right)^m+2M\left(r+(1-r)\ln(1-r)-\sum_{n=2}^{N-1}\frac{r^n}{n(n-1)}\right)\\&\quad+\mu \Phi^N_{M,t}(r)+4M^2\lambda\left(1+\frac{r}{1-r}\right) G_{2,t}(r)+P_q\bigg(\frac{F_{M}(r)}{1-F_{M}(r)}\bigg)\\&\quad-1-2M\left(1-2\ln 2\right).
	\end{align*}
	The constant $ R^{*, \lambda_1,\lambda_2,\cdots\lambda_q}_{\beta,\mu, \lambda,m,N}(M) $ is best possible.
\end{thm}
As a corollary of Theorem \ref{th-3.2}, we obtain the following result in which the cases for $ N=1, 2, 3, 4 $ are discussed completely.
\begin{cor} Let $ f\in \mathcal{P}^{0}_{\mathcal{H}}(M) $ be given by \eqref{e-1.3} and $ 0\leq M <1/(2(\ln 4-1))$, $ \mu,\; \lambda\in\mathbb{R}_{\geq 0} $.
	\begin{enumerate}
		\item[(i)] If $ N=1 $, then $ \mathcal{S}^{f,\lambda_1,\lambda_2,\cdots\lambda_q}_{\beta, \mu, \lambda, m, 1}(r)\leq {d}\left(f(0),\partial \mathbb{D}\right) $
		for $ |z|=r\leq R^{*, \lambda_1,\lambda_2,\cdots\lambda_q}_{\beta,\mu, \lambda,m,1}(M)$, where $ R^{*, \lambda_1,\lambda_2,\cdots\lambda_q}_{\beta,\mu, \lambda,m,1}(M)  $ is the unique root in $ (0,1) $ of the equation 
		\begin{align*}
			r+2M\mathcal{J}_3(r)+\mathcal{J}^M_2(m, r)+\lambda\left(1+\frac{r}{1-r}\right)\mathcal{J}^M_1(r)+P_q\bigg(\frac{F_{M}(r)}{1-F_{M}(r)}\bigg)=0.
		\end{align*}
		\item[(ii)] If $ N=2 $, then $ \mathcal{S}^{f,\lambda_1,\lambda_2,\cdots\lambda_q}_{\beta, \mu, \lambda, m, 2}(r)\leq {d}\left(f(0),\partial \mathbb{D}\right) $
		for $ |z|=r\leq  R^{*, \lambda_1,\lambda_2,\cdots\lambda_q}_{\beta,\mu, \lambda,m,2}(M) $, where $ R^{*, \lambda_1,\lambda_2,\cdots\lambda_q}_{\beta,\mu, \lambda,m,2}(M) $ is the unique root in $ (0,1) $ of the equation 
		\begin{align*}
			2M\mathcal{J}_3(r)+\mathcal{J}^M_2(m, r)+\lambda\left(1+\frac{r}{1-r}\right)\mathcal{J}^M_1(r)+P_q\bigg(\frac{F_{M}(r)}{1-F_{M}(r)}\bigg)=0.
		\end{align*}
		\item[(iii)] If $ N=3 $, then $ \mathcal{S}^{f,\lambda_1,\lambda_2,\cdots\lambda_q}_{\beta, \mu, \lambda, m, 3}(r)\leq {d}\left(f(0),\partial \mathbb{D}\right) $
		for $ |z|=r\leq  R^{*, \lambda_1,\lambda_2,\cdots\lambda_q}_{\beta,\mu, \lambda,m,3}(M)$, where $ R^{*, \lambda_1,\lambda_2,\cdots\lambda_q}_{\beta,\mu, \lambda,m,3}(M) $ is the unique root in $ (0,1) $ of the equation 
		\begin{align*}
			&2M\left(\mathcal{J}_3(r)-\frac{r^2}{2}\right)+\mathcal{J}^M_2(m, r)+\mu\frac{r^3}{1-r}+\lambda\left(1+\frac{r}{1-r}\right)\left(\mathcal{J}^M_1(r)-r^2\right)\\&\quad+P_q\bigg(\frac{F_{M}(r)}{1-F_{M}(r)}\bigg)=0.
		\end{align*}
		\item[(iv)] If $ N=4 $, then $ \mathcal{S}^{f,\lambda_1,\lambda_2,\cdots\lambda_q}_{\beta, \mu, \lambda, m, N}(r)\leq {d}\left(f(0),\partial \mathbb{D}\right) $
		for $ |z|=r\leq  R^{*, \lambda_1,\lambda_2,\cdots\lambda_q}_{\beta,\mu, \lambda,m,4}(M)$, where $ R^{*, \lambda_1,\lambda_2,\cdots\lambda_q}_{\beta,\mu, \lambda,m,4}(M)  $ is the unique root in $ (0,1) $ of the equation 
		\begin{align*}
			&2M\left(\mathcal{J}_3(r)-\frac{r^2}{2}-\frac{r^3}{6}\right)+\mathcal{J}^M_2(m, r)+\mu\frac{r^4}{1-r}+\lambda\left(1+\frac{r}{1-r}\right)\left(\mathcal{J}^M_1(r)-r^2\right)\\&\quad+P_q\bigg(\frac{F_{M}(r)}{1-F_{M}(r)}\bigg)=0.
		\end{align*}
	\end{enumerate}
	The constants $ R^{*, \lambda_1,\lambda_2,\cdots\lambda_q}_{\beta,\mu, \lambda,m,1}(M)$, $R^{*,\lambda_1,\lambda_2,\cdots\lambda_q}_{\beta,\mu, \lambda,m,2}(M)$,  $  R^{*, \lambda_1,\lambda_2,\cdots\lambda_q}_{\beta,\mu, \lambda,m,3}(M)$ and $  R^{*, \lambda_1,\lambda_2,\cdots\lambda_q}_{\beta,\mu, \lambda,m,4}(M) $ all are best possible.
\end{cor}

\begin{rem} 
	For particular choice of \;$\beta=1 $ and $\lambda_j=0 \;\forall \;j=1,2,\cdots q$\;\; and $ \mu,\; \lambda\in\mathbb{R}_{\geq 0}$\;  in Theorem \ref{th-3.1} and \ref{th-3.2}, then Theorem \ref{th-3.1} and  \ref{th-3.2} coincides with \cite[Theorem 2.2]{Aha-CVEE-2022}. Hence  \cite[Theorem 2.2]{Aha-CVEE-2022} becomes a particular case of the Theorem 3.1 and 3.2.
\end{rem}
The following corollary partially addresses the Question \ref{Q2}, while also serving as a harmonic counterpart to Theorem \ref{th-1.18} for the class $\mathcal{P}^{0}_{\mathcal{H}}(M)$.

\begin{cor}	Let $ f\in \mathcal{P}^{0}_{\mathcal{H}}(M) $ be given by \eqref{e-1.3} and $ 0\leq M <1/(2(\ln 4-1))$. Then  we have
	\begin{enumerate}
		\item[(a)] \begin{align*}
			r+\sum_{n=2}^{\infty}\left(|a_n|+|b_n|\right)r^n+\frac{16}{9}\bigg(\frac{S_r}{\pi-S_r}\bigg)\leq {d}\left(f(0),\partial \mathbb{D}\right)
		\end{align*}
		for $ |z|=r\leq R^{*,16/9,0,\cdots0}_{0,\mu, 0,m,1}(M) $, where $R^{*,16/9,0,\cdots0}_{0,\mu, 0,m,1}(M)  $ is the unique root in $ (0,1) $ of the equation 
		\begin{align}\label{EQN-3.25}
			r-1+2M\big((1-r)\big(\log(1-r)-1\big)+\ln4\big)+\dfrac{16}{9}\bigg(\frac{F_{M}(r)}{1-F_{M}(r)}\bigg)=0.
		\end{align}
		The constant $R^{*,16/9,0,\cdots0}_{0,\mu, 0,m,1}(M)$ is best possible. Moreover, we have
		\item[(b)]
		\begin{align*}
			r+\sum_{n=2}^{\infty}\left(|a_n|+|b_n|\right)r^n+\frac{9}{8}\bigg(\frac{S_r}{\pi-S_r}\bigg)\leq {d}\left(f(0),\partial \mathbb{D}\right)
		\end{align*}
		for $ |z|=r\leq R^{*,9/8,0,\cdots0}_{0,\mu, 0,m,1}(M) $, where $R^{*,9/8,0,\cdots0}_{0,\mu, 0,m,1}(M) $ is the unique root in $ (0,1) $ of the equation
		\begin{align}\label{EQN-3.26}
			r-1+2M\big((1-r)\big(\log(1-r)-1\big)+\ln4\big)+\dfrac{9}{8}\bigg(\frac{F_{M}(r)}{1-F_{M}(r)}\bigg)=0.
		\end{align}
		The constant $R^{*,9/8,0,\cdots0}_{0,\mu, 0,m,1}(M)$ is best possible.
	\end{enumerate}
\end{cor}

\begin{figure}
	\begin{center}
		\includegraphics[width=0.5\linewidth]{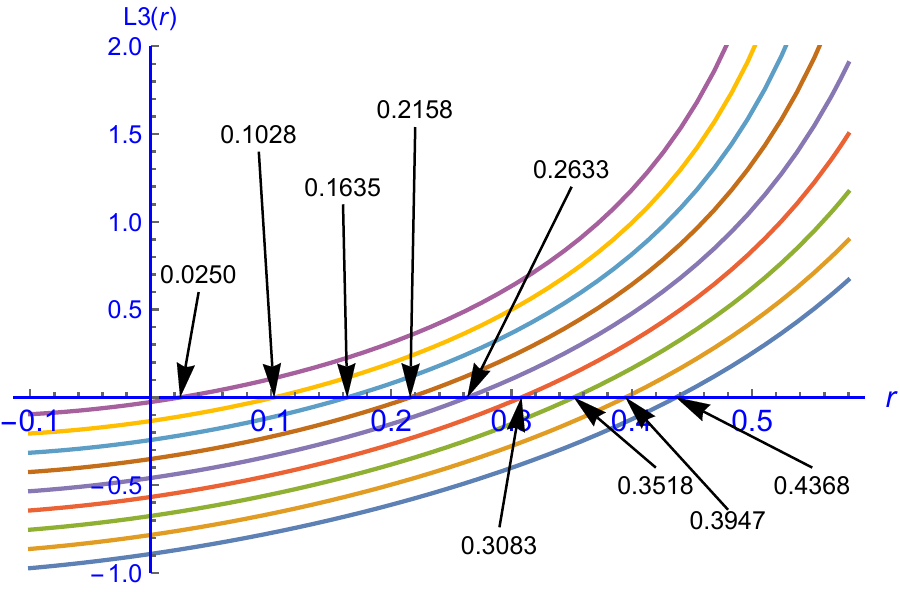}\;\;\;\;\;\;\includegraphics[width=0.5\linewidth]{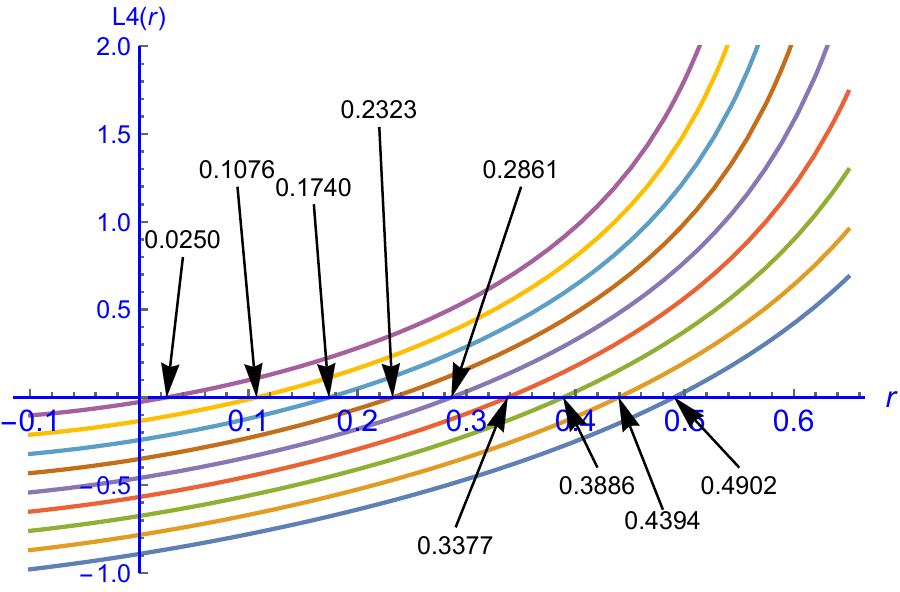}
	\end{center}
	\caption{The figure exhibits the roots of equation \eqref{EQN-3.25}  and the roots of equation \eqref{EQN-3.26}  as defined in Table \ref{T-6} .}
\end{figure}
\begin{table}[ht]
	\centering
	\begin{tabular}{|l|l|l|l|l|l|l|l|l|l|}
		\hline
		$M$& $0.14 $&$0.28 $& $0.42 $& $0.56$& $0.70 $&$0.84 $&$ 0.98 $ &$1.12$& $1.26$  \\
		\hline
		$R^{*,16/9,0,\cdots0}_{0,\mu, 0,m,1}(M)$& $0.4368$&$0.3947 $& $0.3518$& $0.3083$& $0.2633 $&$0.2158$& $0.1635$& $0.1028$ &$ 0.0246 $\\
		\hline
		$ R^{*,9/8,0,\cdots0}_{0,\mu, 0,m,1}(M)$& $0.4902$&$0.4394 $& $0.3886$& $0.3377$& $0.2861 $&$0.2323 $& $0.1740$& $0.1076$ &$ 0.0250 $\\
		\hline
	\end{tabular}\vspace{2.5mm}
	\caption{The table exhibits the roots of the equations \eqref{EQN-3.5} and \eqref{EQN-3.6} for different values of $M.$}
	\label{T-6}
\end{table}

The following corollary of Theorem \ref{th-3.2} is a harmonic analogue of the Theorem \ref{th-2.11} for the class $\mathcal{P}^{0}_{\mathcal{H}}(M).$ Thus a part of Question \ref{Q2} is answered.
\begin{cor}Let $ f\in \mathcal{P}^{0}_{\mathcal{H}}(M) $ be given by \eqref{e-1.3} and $ 0\leq M <1/(2(\ln 4-1))$. Then we have
	\begin{align*}
		|f(z)|+r+\sum_{n=2}^{\infty}\left(|a_n|+|b_n|\right)r^n+2(\sqrt{5}-1)\bigg(\frac{S_r}{\pi-S_r}\bigg)\leq {d}\left(f(0),\partial \mathbb{D}\right)
	\end{align*}  
	for $ |z|=r\leq R^{*,2(\sqrt{5}-1),0,\cdots0}_{1, \mu, 0, 1, 1}(M) $, where $R^{*,2(\sqrt{5}-1),0,\cdots0}_{1, \mu, 0, 1, 1}(M)$ is the unique root  in $ (0,1) $ of the equation 
	\begin{align}\label{eq-3.13}
		r-1+2M\big(r-1+\ln4+(1-r)\ln(1-r)\big)+2(\sqrt{5}-1)\bigg(\frac{F_M(r)}{1-F_M(r)}\bigg)=0.
	\end{align} 
	The constant $R^{*,2(\sqrt{5}-1),0,\cdots0}_{1, \mu, 0, 1, 1}(M)$ is best possible.
\end{cor}
\begin{figure}
	\begin{center}
		\includegraphics[width=0.5\linewidth]{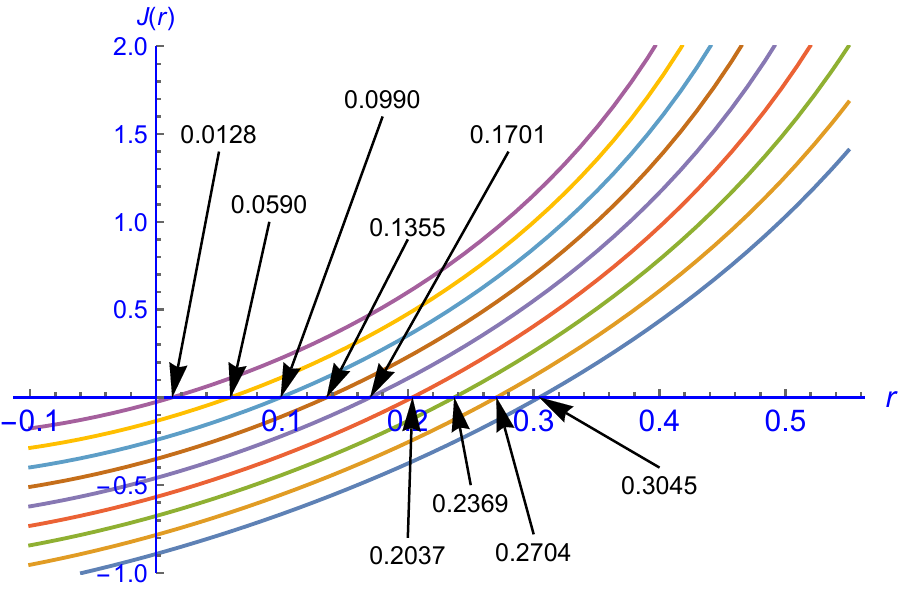}
	\end{center}
	\caption{The figure exhibits the roots of equation \eqref{eq-3.13}    as defined in Table \ref{T-7} .}
\end{figure}
\begin{table}[ht]
	\centering
	\begin{tabular}{|l|l|l|l|l|l|l|l|l|l|}
		\hline
		$M$& $0.14 $&$0.28 $& $0.42 $& $0.56$& $0.70 $&$0.84 $&$ 0.98 $ &$1.12$& $1.26$  \\
		\hline
		$R^{*,2(\sqrt{5}-1),0,\cdots0}_{1, \mu, 0, 1, 1}(M) $& $0.3045$&$0.2704 $& $0.2369$& $0.2037$& $0.1701 $&$0.1355$& $0.0990$& $0.0590$ &$ 0.0128 $\\
		\hline
	\end{tabular}\vspace{2.5mm}
	\caption{The table exhibits the roots of the equations \eqref{eq-3.13}  for different values of $M.$}
	\label{T-7}
\end{table}

\section{Proof of Theorems \ref{th-3.1} and \ref{th-3.2}}

\begin{proof}[\bf Proof of Theorem \ref{th-3.1}]
Let $ f\in \mathcal{P}^{0}_{\mathcal{H}}(M) $ be given by \eqref{e-1.3} and $ 0\leq M <1/(2(\ln 4-1))$. By a straightforward computation, it can be easily shown that
\[
\begin{cases}
	\displaystyle\sum_{n=2}^{\infty}\frac{r^n}{n(n-1)}=r+(1-r)\ln (1-r);\vspace{2mm}\\
	\displaystyle\sum_{n=N}^{\infty}\frac{r^n}{n(n-1)}=r+(1-r)\ln(1-r)-\sum_{n=2}^{N-1}\frac{r^n}{n(n-1)};\vspace{2mm}\\
	\displaystyle \sum_{n=2}^{\infty}\frac{r^{2n}}{n^2(n-1)^2}=\left(r^2+1\right){\rm Li_2\left(r^2\right)}+2\left(r^2-1\right)\ln\left(1-r^2\right)-3r^2;\vspace{2mm}\\
	\displaystyle \sum_{n=t+1}^{\infty}\frac{r^{2n}}{n^2(n-1)^2}=\left(r^2+1\right){\rm Li_2\left(r^2\right)}+2\left(r^2-1\right)\ln\left(1-r^2\right)-3r^2-\sum_{n=2}^{t}\frac{r^{2n}}{n^2(n-1)^2};\vspace{2mm}\\
	\displaystyle
	\sum_{n=N}^{\infty}\frac{2M}{n(n-1)}r^n=\displaystyle2M\left(r+(1-r)\ln(1-r)-\sum_{n=2}^{N-1}\frac{r^n}{n(n-1)}\right).\vspace{2mm}
	
\end{cases}
\]
	In view of growth estimates of Lemma  \ref{lem-2.20} and the estimates above lead to the fact that
	\begin{align}\label{eq-2.8}
		|f(z)|\leq r+2M\sum_{n=2}^{\infty}\frac{r^n}{n(n-1)}=r+2M\left(r+(1-r)\ln (1-r)\right)=: G_{M}(r).
	\end{align}
Since $ N\geq 3 $, then using Lemma \ref{lem-2.19} and the above estimates lead to the following inequality
\begin{align}\label{eq-2.9}
	\sum_{n=N}^{\infty}\left(|a_n|+|b_n|\right)r^n\leq 2M\left(r+(1-r)\ln(1-r)-\sum_{n=2}^{N-1}\frac{r^n}{n(n-1)}\right),
\end{align}
\begin{align}\label{eq-2.10}
	\sum_{n=t+1}^{\infty}\left(|a_n|+|b_n|\right)^2r^{2n}\leq 4M^2\sum_{n=t+1}^{\infty}\frac{r^{2n}}{n^2(n-1)^2}=4M^2 G_{t}(r),
\end{align}
where 
\begin{align*}
	G_{t}(r):=\left(r^2+1\right){\rm Li_2\left(r^2\right)}+2\left(r^2-1\right)\ln\left(1-r^2\right)-3r^2-\sum_{n=2}^{t}\frac{r^{2n}}{n^2(n-1)^2}
\end{align*}and
\begin{align}\label{eq-2.11}
	sgn(t)\sum_{n=1}^{t}\left(|a_n|+|b_n|\right)^2\frac{r^N}{1-r}\leq sgn(t)\frac{4M^2r^N}{1-r}\sum_{n=1}^{t}\frac{1}{n^2(n-1)^2}=:\Phi^N_{M,t}(r).
\end{align}
We need to find upper bound of the quantity $S_r/\pi$ for the class $\mathcal{P}^0_H(M).$ It is well known that (see \cite[p. 113]{Duren-2004}) the area of the disk $\mathbb{D}_r:=\{z\in\mathbb{C}:|z|<r\}$ under the harmonic map $ f=h+\bar{g} $ is given by
\begin{align}\label{e-3.9}
S_r=\iint\limits_{\mathbb{D}_r}J_f(z)dxdy=\iint\limits_{\mathbb{D}_r}\left(|{h^{\prime}(z)}|^2-|{g^{\prime}(z)}|^2\right)dxdy.
\end{align}
Therefore, in view of Lemma \ref{lem-2.19}, it follows that (see \cite{Ahamed-AMP-2021,Ahamed-CVEE-2021} for detailed information)
\begin{align}\label{eq-2.3}	\frac{S_r}{\pi}&=r^2+\sum_{n=2}^{\infty}n\left(|a_n|^2-|b_n|^2\right)r^{2n}\\&\nonumber\leq r^2+4M^2\sum_{n=2}^{\infty}\frac{r^{2n}}{n(n-1)^2}\\&= r^2+4M[r^2Li_2(r^2)-(r^2+(1-r^2)\log(1-r^2)):=F_{M}(r)\nonumber.
\end{align}
By a simple computation, the following can be shown
\begin{align}\label{eq-3.19}
	\dfrac{d}{dr}\big(F_{M}(r)\big)>0 \; \mbox{for} \; r\in(0,1).
\end{align}
Thus, using \eqref{eq-2.8} to \eqref{eq-2.11} and \eqref{eq-2.3}, a straight forward computation yeilds that 
\begin{align} \label{eq-2.13}
	&G^{f,\lambda_1,\lambda_2,\cdots\lambda_q}_{\mu, \lambda, m, N}(r)\\\nonumber&\leq \ \beta\left(G_{M}(r)\right)^m+2M\left(r+(1-r)\ln(1-r)-\sum_{n=2}^{N-1}\frac{r^n}{n(n-1)}\right)+\mu \Phi^N_{M,t}(r)\\&\nonumber\quad+4M^2\lambda\left(1+\frac{r}{1-r}\right) G_{t}(r)+P_q(F_{M}(r))\\&\nonumber\leq 1+2M\left(1-2\ln 2\right)
\end{align} 
for $ |z|=r\leq R^{\lambda_1,\lambda_2,\cdots\lambda_q}_{\beta,\mu, \lambda,m,N}(M)$, where $ R^{\lambda_1,\lambda_2,\cdots\lambda_q}_{\beta,\mu, \lambda,m,N}(M) $ is the smallest root in $ (0, 1) $ of equation $ \Phi^{\lambda_1,\lambda_2,\cdots\lambda_q}_{\beta,\mu, \lambda,m,N}(r)=0 $ .  It is not hard to show that $ R^{\lambda_1,\lambda_2,\cdots\lambda_q}_{\beta,\mu, \lambda,m,N}(M)$ is the unique root of the equation $ \Phi^{\lambda_1,\lambda_2,\cdots\lambda_q}_{\beta,\mu, \lambda,m,N}(r)=0 $ in $ (0, 1) $, where $\Phi^{\lambda_1,\lambda_2,\cdots\lambda_q}_{\beta,\mu, \lambda,m,N}(r)$ is real-valued differentiable function on $(0,1)$ is defined by \eqref{eq-2.4}. Indeed, a simple computation for all $r\in(0,1)$ shows that 
\[
\begin{cases}
	\displaystyle\frac{d}{dr}\left(G_{M}(r)\right)^m=m\bigg(r+2M\sum_{n=2}^{\infty}\frac{r^n}{n(n-1)}\bigg)^{m-1}\bigg(1+2M\sum_{n=2}^{\infty}\frac{nr^{n-1}}{n(n-1)}\bigg)>0;\vspace{2mm}\\
	\displaystyle \displaystyle\frac{d}{dr}\bigg(2M\bigg(r+(1-r)\ln(1-r)-\sum_{n=2}^{N-1}\frac{r^n}{n(n-1)}\bigg)\bigg)=\sum_{n=N}^{\infty}\frac{2Mnr^{n-1}}{n(n-1)}>0;\vspace{2mm}\\
	\displaystyle \frac{d}{dr}(\Phi^N_{M,t}(r))=\displaystyle \frac{4M^2r^{N-1}(n(1-r)+r)}{(1-r)^2}sgn(t)\sum_{n=1}^{t}\frac{1}{n^2(n-1)^2}>0;\vspace{2mm}\\
  \displaystyle\frac{d}{dr}\bigg(\left(1+\frac{r}{1-r}\right) G_{t}(r)\bigg)\\=\bigg(1+\dfrac{1}{1-r}\bigg)\displaystyle\sum_{n=t+1}^{\infty}\frac{{2n}r^{2n-1}}{n^2(n-1)^2}+\displaystyle\sum_{n=t+1}^{\infty}\frac{r^{2n}}{n^2(n-1)^2}\frac{1}{(1-r)^2}>0.

\end{cases}
\]
In view of \eqref{eq-3.19} and above equations, it is easy to see that
\begin{align*}
	&\frac{d}{dr}\left(\Phi^{\lambda_1,\lambda_2,\cdots\lambda_q}_{\beta,\mu,\lambda,m,N}(r)\right)\\=\nonumber&\beta\frac{d}{dr}\big(\left(G_{M}(r)\right)^m\big)+\dfrac{d}{dr}\bigg(2M\left(r+(1-r)\ln(1-r)-\sum_{n=2}^{N-1}\frac{r^n}{n(n-1)}\right)\bigg)\\&\quad+\mu\dfrac{d}{dr}\big(\Phi^N_{M,t}(r)\big)+4M^2\lambda\;\dfrac{d}{dr}\bigg(\left(1+\frac{r}{1-r}\right) G_{t}(r)\bigg)\nonumber\nonumber+\lambda_1\frac{d}{dr}(F_{M}(r))\\&\quad+2\lambda_2F_{M}(r)\frac{d}{dr}(F_{M}(r))+\cdots+q\lambda_q(F_{M}(r))^{q-1}\times\nonumber\nonumber\frac{d}{dr}(F_{M}(r))>0\nonumber \;\mbox{for all}\;r\in (0, 1).
\end{align*}
Evidently, $\Phi^{\lambda_1,\lambda_2,\cdots\lambda_q}_{\beta,\mu,\lambda,m,N}(r)$ is an increasing function of $r$ on $(0,1)$. Moreover,  $\Phi^{\lambda_1,\lambda_2,\cdots\lambda_q}_{\beta,\mu,\lambda,m,N}(0)= -1-2M\left(1-2\ln 2\right)<0$ and $\lim\limits_{{r\rightarrow1}}\Phi^{\lambda_1,\lambda_2,\cdots\lambda_q}_{\beta,\mu,\lambda,m,N}(r)=+\infty$  and hence the equation $\Phi^{\lambda_1,\lambda_2,\cdots\lambda_q}_{\beta,\mu,\lambda,m,N}(r)=0$  has the unique positive root $R^{\lambda_1,\lambda_2,\cdots\lambda_q}_{\beta,\mu, \lambda,m,N}(M)$ in $(0,1).$
 Therefore, we have
\begin{align}\label{eq-2.16} \Phi^{\lambda_1,\lambda_2,\cdots\lambda_q}_{\beta,\mu,\lambda,m,N}\big(R^{\lambda_1,\lambda_2,\cdots\lambda_q}_{\beta,\mu, \lambda,m,N}(M)\big)+1+2M\left(1-2\ln 2\right)=1+2M\left(1-2\ln 2\right).
\end{align}
The Euclidean distance between $f(0)$ and the boundary of $f(\mathbb{D})$ is given by 
\begin{equation} \label{e-2.1}
	d(f(0), \partial f(\mathbb{D}))= \liminf \limits_{|z|=r\rightarrow 1} |f(z)-f(0)|.
\end{equation}
Since $f(0)=0$, in view of Lemma \ref{lem-2.20} and \eqref{e-2.1}, a simple computation shows that
\begin{equation} \label{eq-2.17}
	d(f(0), \partial f(\mathbb{D}))=\liminf \limits_{|z|=r\rightarrow 1} |f(z)|\geq 1+2M\sum\limits_{n=2}^{\infty}  \frac{(-1)^{n-1}}{n(n-1)}=1+2M\left(1-2\ln 2\right).
\end{equation}

Therefore, in view of \eqref{eq-2.13} and \eqref{eq-2.17}, we see that the desired inequality $ G^{f,\lambda_1,\lambda_2,\cdots\lambda_q}_{\beta, \mu, \lambda, m, N}(r) \leq {d}\left(f(0),\partial \mathbb{D}\right) $ holds for $r\leq R^{\lambda_1,\lambda_2,\cdots\lambda_q}_{\beta,\mu, \lambda,m,N}(M)$.\vspace{1.2mm}

Next part of the proof is show that the radius $R^{\lambda_1,\lambda_2,\cdots\lambda_q}_{\beta,\mu, \lambda,m,N}(M)$ is best possible. In order to that, we consider the function $ f=f_{M} $ defined by 
\begin{equation}\label{Eq-2.16}
	f_{M}(z)=z+\sum_{n=2}^{\infty}\frac{2Mz^n}{n(n-1)}.
\end{equation}
We see that $ f_{M}\in\mathcal{P}^{0}_{\mathcal{H}}(M) $ and  for $ f=f_{M} $, in view of \eqref{e-2.1}, by a routine computation, we obtain the explicit form of the Euclidean distance as
\begin{equation}\label{e-2.2424}
	d(f_{M}(0),\partial f_{M}(\mathbb{D}))=1+2M\left(1-2\ln 2\right).
\end{equation}
Therefore, for $ f=f_{M} $ and $z=r>R^{\lambda_1,\lambda_2,\cdots\lambda_q}_{\beta,\mu, \lambda,m,N}(M)$,  in view of \eqref{eq-2.16} 
and \eqref{e-2.2424}, it can be easily shown that 
\begin{align*} 
	&G^{f_M,\lambda_1,\lambda_2,\cdots\lambda_q}_{\mu, \lambda, m, N}(r)\\&=\beta\bigg(r+2M\sum_{n=2}^{\infty}\frac{r^n}{n(n-1)}\bigg)^m+\sum_{n=N}^{\infty}\frac {2M}{n(n-1)}r^n+\mu\; sgn(t)\sum_{n=1}^{t}\frac {4M^2}{n^2(n-1)^2}\frac{r^N}{1-r}\\&\quad+\lambda\left(1+\frac{r}{1-r}\right)\sum_{n=t+1}^{\infty}\frac {4M^2}{n^2(n-1)^2}r^{2n}+P_q\bigg(r^2+4M^2\sum_{n=2}^{\infty}\frac{r^{2n}}{n(n-1)^2}\bigg)\\&>\beta\bigg(R^{\lambda_1,\lambda_2,\cdots\lambda_q}_{\beta,\mu, \lambda,m,N}(M)+2M\sum_{n=2}^{\infty}\frac{\big(R^{\lambda_1,\lambda_2,\cdots\lambda_q}_{\beta,\mu, \lambda,m,N}(M)\big)^n}{n(n-1)}\bigg)^m+\sum_{n=N}^{\infty}\frac {2M}{n(n-1)}\big(R^{\lambda_1,\lambda_2,\cdots\lambda_q}_{\beta,\mu, \lambda,m,N}(M)\big)^n\\&\quad+\mu\; sgn(t)\sum_{n=1}^{t}\frac {4M^2}{n^2(n-1)^2}\frac{\big(R^{\lambda_1,\lambda_2,\cdots\lambda_q}_{\beta,\mu, \lambda,m,N}(M)\big)^N}{1-R^{\lambda_1,\lambda_2,\cdots\lambda_q}_{\beta,\mu, \lambda,m,N}(M)}+\lambda\left(1+\frac{R^{\lambda_1,\lambda_2,\cdots\lambda_q}_{\beta,\mu, \lambda,m,N}(M)}{1-R^{\lambda_1,\lambda_2,\cdots\lambda_q}_{\beta,\mu, \lambda,m,N}(M)}\right)\times\\&\quad\sum_{n=t+1}^{\infty}\frac {4M^2}{n^2(n-1)^2}\big(R^{\lambda_1,\lambda_2,\cdots\lambda_q}_{\beta,\mu, \lambda,m,N}(M)\big)^{2n}+P_q\bigg(\big(R^{\lambda_1,\lambda_2,\cdots\lambda_q}_{\beta,\mu, \lambda,m,N}(M)\big)^2+4M^2\sum_{n=2}^{\infty}\frac{\big(R^{\lambda_1,\lambda_2,\cdots\lambda_q}_{\beta,\mu, \lambda,m,N}(M)\big)^{2n}}{n(n-1)^2}\bigg)\\&=\Phi^{\lambda_1,\lambda_2,\cdots\lambda_q}_{\beta,\mu,\lambda,m,N}\big(R^{\lambda_1,\lambda_2,\cdots\lambda_q}_{\beta,\mu, \lambda,m,N}(M)\big)+1+2M\left(1-2\ln 2\right)\\&=1+2M\left(1-2\ln 2\right)\\&=d(f_{M}(0),\partial f_{M}(\mathbb{D})).
\end{align*}
Clearly, the constant $R^{\lambda_1,\lambda_2,\cdots\lambda_q}_{\beta,\mu, \lambda,m,N}(M) $ is best possible. This completes the proof.
\end{proof}	

\begin{proof}[\bf Proof of Theorem \ref{th-3.2}]
	Let $ f\in \mathcal{P}^{0}_{\mathcal{H}}(M) $ be given by \eqref{e-1.3} and $ 0\leq M <1/(2(\ln 4-1))$. Since $S_r/\pi\leq F_M(r)$, an easy computation shows that
	\begin{align*}
		\frac{S_r}{\pi-S_r}\leq \frac{F_{M}(r)}{1-F_{M}(r)}\quad\mbox{for}\; r\in (0, 1).
	\end{align*}	
	From the immediate last inequality and \eqref{eq-2.8} to \eqref{eq-2.11}, a straightforward computation shows that
	\begin{align} \label{EQ-3.28}
		&\mathcal{S}^{f,\lambda_1,\lambda_2,\cdots\lambda_q}_{\mu, \lambda, m, N}(r)\\\nonumber&\leq \beta\left(r+2M\sum_{n=2}^{\infty}\frac{r^n}{n(n-1)}\right)^m+2M\sum_{n=N}^{\infty}\frac{r^n}{n(n-1)}+\mu \;sgn(t)\frac{4M^2r^N}{1-r}\sum_{n=1}^{t}\frac{1}{n^2(n-1)^2}\\&\nonumber\quad+\lambda\left(1+\frac{r}{1-r}\right)\sum_{n=t+1}^{\infty}\frac{4M^2r^{2n}}{n^2(n-1)^2}+P_q\bigg(\dfrac{F_{M}(r)}{1-F_{M}(r)}\bigg)\\&\nonumber\leq  \beta\left(G_{M}(r)\right)^m+2M\left(r+(1-r)\ln(1-r)-\sum_{n=2}^{N-1}\frac{r^n}{n(n-1)}\right)+\mu \Phi^N_{M,t}(r)\\&\nonumber\quad+4M^2\lambda\left(1+\frac{r}{1-r}\right) G_{t}(r)+P_q\bigg(\dfrac{F_{M}(r)}{1-F_{M}(r)}\bigg)\\&\nonumber\leq 1+2M\left(1-2\ln 2\right)
	\end{align} 
	for $ |z|=r\leq R^{*, \lambda_1,\lambda_2,\cdots\lambda_q}_{\beta,\mu, \lambda,m,N}(M) $, where $ R^{*, \lambda_1,\lambda_2,\cdots\lambda_q}_{\beta,\mu, \lambda,m,N}(M) $ is the smallest root of $ \Phi^{*, \lambda_1,\lambda_2,\cdots\lambda_q}_{\beta,\mu, \lambda,m,N}(r)=0 $ in $ (0, 1) $. As
	 \begin{align}\label{Eq-1.9}
		\dfrac{d}{dr}(F_M(r))>0,\;\mbox{then}\;\;\frac{d}{dr}\bigg(\frac{G(r)}{1-G(r)}\bigg)=\displaystyle \dfrac{\frac{d}{dr}\big(G(r)\big)}{(1-G(r))^2}>0.
	\end{align}  
	By the similar argument as in proof of Theorem \ref{th-3.1}, it can be shown that $ R^{*, \lambda_1,\lambda_2,\cdots\lambda_q}_{\beta,\mu, \lambda,m,N}(M) $ is the unique root of the equation $ \Phi^{*, \lambda_1,\lambda_2,\cdots\lambda_q}_{\beta,\mu, \lambda,m,N}(r)=0 $ in $ (0, 1) $. Therefore, we have
	\begin{align}\label{e-2.22}&  \Phi^{*, \lambda_1,\lambda_2,\cdots\lambda_q}_{\beta,\mu, \lambda,m,N}\bigg(R^{*, \lambda_1,\lambda_2,\cdots\lambda_q}_{\beta,\mu, \lambda,m,N}(M)\bigg)1+2M\left(1-2\ln 2\right)=1+2M\left(1-2\ln 2\right).
	\end{align}
	
	Therefore, in view of  \eqref{eq-2.17} and \eqref{EQ-3.28} , we obtain the desired inequality $ \mathcal{S}^{f,\lambda_1,\lambda_2,\cdots\lambda_q}_{\beta, \mu, \lambda, m, N}(r) \leq {d}\left(f(0),\partial \mathbb{D}\right) $\; for $r\leq R^{*, \lambda_1,\lambda_2,\cdots\lambda_q}_{\beta,\mu, \lambda,m,N}(M)$.\vspace{1.2mm} 
	
	In order to complete the proof, it suffices to show that the constant $ R^{*, \lambda_1,\lambda_2,\cdots\lambda_q}_{\beta,\mu, \lambda,m,N}(M) $ is best possible. Henceforth, we consider the function defined by \eqref{Eq-2.16}.
	For $ f=f_{M} $ and $z=r>R^{*, \lambda_1,\lambda_2,\cdots\lambda_q}_{\beta,\mu, \lambda,m,N}(M) $, by the similar argument as in the proof of Theorem \ref{th-3.1}, using \eqref{e-2.2424} 
	and \eqref{e-2.22}, it can be shown that 
	\begin{align*} 
		&\mathcal{S}^{f_M,\lambda_1,\lambda_2,\cdots\lambda_q}_{\mu, \lambda, m, N}(r)\\&>\ \Phi^{*, \lambda_1,\lambda_2,\cdots\lambda_q}_{\beta,\mu, \lambda,m,N}\bigg(R^{*, \lambda_1,\lambda_2,\cdots\lambda_q}_{\beta,\mu, \lambda,m,N}(M)\bigg)1+2M\left(1-2\ln 2\right)=1+2M\left(1-2\ln 2\right)\\&=1+2M\left(1-2\ln 2\right)\\&=d(f_{M}(0),\partial f_{M}(\mathbb{D})).
	\end{align*}
	Therefore, $ R^{*, \lambda_1,\lambda_2,\cdots\lambda_q}_{\beta,\mu, \lambda,m,N}(M)$ is best possible. This completes the proof.
\end{proof}	
\noindent{\bf Acknowledgment:} The authors would like to thank the anonymous referees for their helpful suggestions and comments to enhance the clarity and presentation of the paper. \\
\vspace{1.2mm}

\noindent\textbf{Compliance of Ethical Standards:}\\

\noindent\textbf{Conflict of interest.} The authors declare that there is no conflict  of interest regarding the publication of this paper.\vspace{1.5mm}

\noindent\textbf{Funds.} No funds.\vspace{1.5mm}

\noindent\textbf{Data availability statement.}  Data sharing not applicable to this article as no datasets were generated or analysed during the current study.

\end{document}